\documentclass[12pt]{article}
\usepackage{amssymb,amsthm,amsmath,amsfonts,latexsym,tikz,hyperref}
\usepackage[hmargin=1in,vmargin=1in]{geometry}

\newtheorem{thm}{Theorem}[section]
\newtheorem{prop}[thm]{Proposition}
\newtheorem{cor}[thm]{Corollary}
\newtheorem{lem}[thm]{Lemma}
\newtheorem{conj}[thm]{Conjecture}
\newtheorem{exa}[thm]{Example}
\newtheorem{question}[thm]{Question}
\newtheorem{problem}[thm]{Problem}

\newcommand{\etalchar}[1]{$^{#1}$}

\DeclareMathOperator{\Asc}{Asc}
\DeclareMathOperator{\dasc}{dasc}
\DeclareMathOperator{\dAsc}{dAsc}
\DeclareMathOperator{\wasc}{wasc}
\DeclareMathOperator{\wAsc}{wAsc}

\DeclareMathOperator{\dMtx}{dMtx}
\DeclareMathOperator{\Mtx}{Mtx}

\DeclareMathOperator{\A}{A}
\DeclareMathOperator{\wA}{wA}
\DeclareMathOperator{\dA}{dA}
\DeclareMathOperator{\tAsc}{tAsc}

\DeclareMathOperator{\act}{act}

\DeclareMathOperator{\pe}{pe}
\DeclareMathOperator{\mx}{mx}
\DeclareMathOperator{\po}{po}
\DeclareMathOperator{\mh}{mh}
\DeclareMathOperator{\rmax}{rmax}
\DeclareMathOperator{\rmin}{rmin}
\DeclareMathOperator{\lef}{lef}
\DeclareMathOperator{\ddis}{ddis}

\DeclareMathOperator{\RDT}{RDT}
\DeclareMathOperator{\dRDT}{dRDT}
\DeclareMathOperator{\I}{I}
\DeclareMathOperator{\M}{M}
\DeclareMathOperator{\dMch}{dMch}
\DeclareMathOperator{\dI}{dI}
\DeclareMathOperator{\di}{di}
\DeclareMathOperator{\tI}{tI}

\DeclareMathOperator{\RGF}{RGF}
\DeclareMathOperator{\rp}{rp}

\newcommand{\bfig}{\begin{figure}}
\newcommand{\efig}{\end{figure}}

\newcommand{\ben}{\begin{enumerate}}
\newcommand{\een}{\end{enumerate}}
\newcommand{\ble}{\begin{lem}}
\newcommand{\ele}{\end{lem}}
\newcommand{\bth}{\begin{thm}}
\renewcommand{\eth}{\end{thm}}
\newcommand{\bpr}{\begin{prop}}
\newcommand{\epr}{\end{prop}}
\newcommand{\bco}{\begin{cor}}
\newcommand{\eco}{\end{cor}}
\newcommand{\bcon}{\begin{conj}}
\newcommand{\econ}{\end{conj}}
\newcommand{\bde}{\begin{defn}}
\newcommand{\ede}{\end{defn}}
\newcommand{\bex}{\begin{exa}}
\newcommand{\eex}{\end{exa}}
\newcommand{\barr}{\begin{array}}
\newcommand{\earr}{\end{array}}
\newcommand{\btab}{\begin{tabular}}
\newcommand{\etab}{\end{tabular}}
\newcommand{\beq}{\begin{equation}}
\newcommand{\eeq}{\end{equation}}
\newcommand{\bea}{\begin{eqnarray*}}
\newcommand{\eea}{\end{eqnarray*}}
\newcommand{\bal}{\begin{align*}}
\newcommand{\bce}{\begin{center}}
\newcommand{\ece}{\end{center}}
\newcommand{\bpi}{\begin{picture}}
\newcommand{\epi}{\end{picture}}
\newcommand{\bpp}{\begin{picture}}
\newcommand{\epp}{\end{picture}}
\newcommand{\bfi}{\begin{figure} \begin{center}}
\newcommand{\efi}{\end{center} \end{figure}}
\newcommand{\bprf}{\begin{proof}}
\newcommand{\eprf}{\end{proof}\medskip}

\newcommand{\bsl}{\begin{slide}{}}
\newcommand{\esl}{\end{slide}}
\newcommand{\bfr}{\begin{frame}}
\newcommand{\efr}{\end{frame}}

\newcommand{\hqed}{\hfill \qed}

\newcommand{\eqed}[1]{$\textcolor{white}{\qed}\hfill{\dil#1}\hfill\qed$}

\newcommand{\ol}{\overline}

\newcommand{\hs}[1]{\hspace{#1}}
\newcommand{\hso}[1]{\hspace{-1pt}}
\newcommand{\vs}[1]{\vspace{#1}}

\newcommand{\emp}{\emptyset}

\newcommand{\sbe}{\subseteq}

\newcommand{\spe}{\supseteq}
\newcommand{\setm}{\setminus}

\newcommand{\case}[4]{\left\{\barr{ll}#1&\mbox{#2}\\#3&\mbox{#4}\earr\right.}

\def\<{\langle}
\def\>{\rangle}

\newcommand{\ra}{\rightarrow}

\newcommand{\al}{\alpha}
\newcommand{\be}{\beta}

\newcommand{\de}{\delta}

\newcommand{\io}{\iota}
\newcommand{\ka}{\kappa}

\newcommand{\si}{\sigma}

\newcommand{\Si}{\Sigma}

\newcommand{\7}{{\bf 7}}

\newcommand{\cS}{{\cal S}}

\newcommand{\fS}{{\mathfrak S}}

\DeclareMathOperator{\Av}{Av}

\DeclareMathOperator{\asc}{asc}

\DeclareMathOperator{\sgn}{sgn}

\newcommand{\dil}{\displaystyle}

\begin{document}
\pagestyle{plain}

\title{Difference ascent sequences
}
\author{Mark Dukes\\[-5pt]
\small School of Mathematics and Statistics, University College Dublin,\\[-5pt]
\small Belfield, Dublin 4, Ireland, {\tt mark.dukes@ucd.ie}\\
Bruce E. Sagan\\[-5pt]
\small Department of Mathematics, Michigan State University,\\[-5pt]
\small East Lansing, MI 48824-1027, USA, {\tt sagan@math.msu.edu}
}

\date{\today\\[10pt]
	\begin{flushleft}
	                                       \\[5pt]
	\end{flushleft}}

\maketitle

\begin{abstract}
Let $\al=a_1 a_2\ldots a_n$ be a sequence of nonnegative integers.
The ascent set  of $\al$, $\Asc\al$,  consists of all indices $k$ where $a_{k+1}>a_k$.
An ascent sequence is $\al$ where the growth of the $a_k$ is bounded by the elements of $\Asc\al$.
These sequences were introduced by Bousquet-M\'elou, Claesson, Dukes and Kitaev and have many wonderful properties.  In particular, they are in bijection with unlabeled $(2+2)$-free posets,
permutations avoiding a particular bivincular pattern, certain upper-triangular nonnegative integer matrices, and a class of matchings.
A weak ascent of $\al$ is an index $k$ with $a_{k+1}\ge a_k$ and weak ascent sequences are defined analogously to ascent sequences.
These were studied by B\'enyi, Claesson and Dukes and shown to have analogous equinumerous sets.
Given a nonnegative integer $d$, we define a difference $d$ ascent to be an index $k$ such 
that $a_{k+1}>a_k-d$.  We study the properties of the corresponding $d$-ascent sequences, showing that some of the maps from the weak case can be extended to bijections for general $d$ while the extensions of others continue to be injective (but not surjective).  We also make connections with other combinatorial objects such as rooted duplication trees and restricted growth functions.
\end{abstract}

\section{Introduction}

For $m,n$  nonnegative integers we let $[m,n]=\{m,m+1,\ldots,n\}$, and when $m=1$ we abbreviate this to $[n]=\{1,2,\ldots,n\}$.
Note that $[n]=\emp$ if $n=0$.
Consider a sequence  $\al=a_1 a_2\ldots a_n$ of nonnegative integers.  
We will sometimes put commas between the $a_k$ for readability.
Also note our convention that we will use Greek letters for sequences and the corresponding roman letters for their elements.
For $k\in[n]$ we will use
$$
\al_k = a_1 a_2 \ldots a_k
$$
for the prefix of $\al$ with $k$ elements and also let $\al_0=\emp$ be the empty sequence.
The {\em ascent set} and {\em ascent number} of $\al$ are
$$
\Asc\al=\{k\in[n-1] \mid a_{k+1} > a_k\}
$$
and
$$
\asc\al=\#\Asc\al,
$$
respectively, where we will use a hash symbol or a pair of vertical bars to denote cardinality.  For example, if
\begin{equation}
\label{al}  
\al = 0,1,1,0,2,1,2,4
\end{equation}
then $\Asc\al=\{1,4,6,7\}$ and $\asc\al=4$.
Call $\al$ an {\em ascent sequence} if 
\begin{enumerate}
    \item[(a1)] $a_1=0$, and
    \item[(a2)] $a_{k+1}\le 1 + \asc\al_k$ for $k\in[n-1]$.
\end{enumerate}
It is easy to check that $\al$ in the previous example is an ascent sequence while
\begin{equation}
\label{al'}
  \al' = 0,1,1,0,3,1,2, 4  
\end{equation}
is not because $\al'_5=3$ while $1+\asc(0,1,1,0)= 2$ which is smaller.
For $n\ge0$ we  will use the notation
$$
\A_n =\{\al=a_1 a_2\ldots a_n \mid \text{$\al$ is an ascent sequence}\}.
$$
For example
$$
\A_3 = \{000,001,010,011,012\}.
$$
Bousquet-M\'elou, Claesson, Dukes and Kitaev~\cite{BMCDK:2fp} were the first to define and study ascent sequences and they have since been considered by numerous authors such 
as~\cite{CDDDS:0aa,CCEG:paa,CCEG:paa,CDK:des,CL:n!m,duk:gbs,DM:rba,DP:asu,DS:paa,FJLYZ:nda,JS:pbs,LF:1ai,MS:err,pud:asb,yan:bis}.
Ascent sequences are known to be in bijection with unlabeled $(2+2)$-free posets,
permutations avoiding a  bivincular pattern of length $3$, certain upper-triangular nonnegative integer matrices, and a class of matchings.

Very recently, B\'enyi, Claesson, and Dukes~\cite{BCD:was} introduced weak ascent sequences.  The set of {\em weak ascents} of $\al$ is
$$
\wAsc\al=\{k\in[n-1] \mid a_{k+1} \ge a_k\}
$$
with corresponding {\em weak ascent number}
$$
\wasc\al=\#\wAsc\al.
$$
Note that $\wAsc\al\spe\Asc\al$.
Taking $\al'$ as in~\eqref{al'} we see that 
$\wAsc\al'=\{1,2,4,6,7\}$.
The sequence $\al$ is a {\em weak ascent sequence} if
 \begin{enumerate}
    \item[(w1)] $a_1=0$, and
    \item[(w2)] $a_{k+1}\le 1 + \wasc\al_k$ for $k\in[n-1]$.
\end{enumerate} 
Even though $\al'$ as above was not an ascent sequence, it is a weak ascent sequence.
In fact, if we let
$$
\wA_n = \{\al=a_1 a_2\ldots a_n \mid \text{$\al$ is a weak ascent sequence}\}.
$$
then clearly $\wA_n\spe\A_n$ for all $n\ge0$.  As an example,
$$
\wA_3=\{000,001,002,010,011,012\}.
$$
In~\cite{BCD:was} the authors showed that weak ascent sequences are equinumerous with factorial posets avoiding a specially labeled $3+1$, permutations avoiding a bivincular pattern of length $4$, upper-triangular $0$-$1$ matrices satisfying a column restriction, and matchings with a restriction on their nestings.

The purpose of the present work is to introduce and study a more general class of sequences which includes both ascent sequences and weak ascent sequences as special cases.  For $d\ge0$ the {\em difference $d$ ascent set} or simply {\em $d$-ascent set} of $\al$ is
$$
\dAsc\al = \{k\in[n-1] \mid a_{k+1} > a_k -d\}.
$$
with corresponding {\em $d$-ascent number}
$$
\dasc\al =\#\dAsc\al.
$$
Note that when $d=0$ and $d=1$ we have $\dAsc\al = \Asc\al$ and $\dAsc\al=\wAsc\al$, respectively.  For example, when $d=2$ we have
$$
\dAsc(0,1,0,2,3,1,4,1,7) = \{1,2,3,4,5,6,8\}.
$$
Unsurprisingly, the definition of a {\em $d$-ascent sequence} $\al$ is one where
 \begin{enumerate}
    \item[(d1)] $a_1=0$, and
    \item[(d2)] $a_{k+1}\le 1 + \dasc\al_k$ for $k\in[n-1]$.
\end{enumerate} 
The example sequence just given is a $2$-ascent sequence.
Let
$$
\dA_n = \{\al=a_1 a_2\ldots a_n \mid \text{$\al$ is a $d$-ascent sequence}\}.
$$
These difference $d$ ascent sequences are different from the $p$-ascent sequences that were introduced by Kiteav and Remmel~\cite{KR2}.
In what follows, we will continue the convention initiated above using notation starting with  $d$ for notions concerning general $d$-ascent sequences,
with $w$ for the corresponding ideas applied to weak ascent sequences, and with no added initial letter for the ascent sequence case.
Sometimes the $d=2$ sequences will be of special interest and then the prefix will be t as in $\tAsc\al$.

The rest of this paper is structured as follows.  In the next section we will give a bijection between $d$-ascent sequence and  upper triangular $0$-$1$ matrices satisfying  a condition on the columns.
This is a generalization of the map given in B\'enyi et al.~\cite[Sec. 3]{BCD:was} for the case of weak ascent sequences.
Section~\ref{matc} is devoted to a bijection between $\dA_n$ and a set of matchings with restricted nestings.
In Section~\ref{per} we show that there is an injection from $d$-ascent sequence to permutations avoiding a bivincular pattern of length $d+3$. 
Next, we construct another injection with domain $\dA_n$ for $d\ge1$, this one to the factorial posets of Claesson and Linusson~\cite{CL:n!m} avoiding a specially labeled poset $P_{d+3}$ with $d+3$ elements.  
Call a $d$-ascent sequence {\em $d$-increasing} if every index (except the last) is a $d$-ascent.  In Section~\ref{did} we show that there is an alternating sum recursion of such sequences.  Also, for $d=2$, they are in bijection with the rooted duplication trees of Gasceul, Hendy, Jean-Marie, and McLachlan~\cite{GHJM:ctd}.
The special factorial posets considered earlier were for $d\ge1$, but one can also define an analogous $P_3$.
We show in Section~\ref{rgf} that the factorial posets avoiding $P_3$ are in bijection with restricted growth functions.  We end with a section of comments and suggestions for future work.


\section{Matrices}
\label{matr}

We will now give a bijection between difference $d$ ascent sequences and matrices satisfying a certain column condition.  Our map is a simplification and generalization of the one for weak ascent sequences given in~\cite{BCD:was}.

\begin{figure}
    \centering
$$
\barr{rc}
M=
&
\left[
\barr{ccccc}
1&0&1&1&1\\
0&1&0&0&0\\
0&0&0&0&1\\
0&0&0&1&0\\
0&0&0&0&1
\earr
\right]
\\
\\
\barr{r}
j:\\
\rmin c_j:\\
\rmax c_j:
\earr
&
\barr{ccccc}
0 & 1 & 2 & 3 & 4\\
0 & 1 & 0 & 0 & 0\\
0 & 1 & 0 & 3 & 4
\earr
\earr
$$
    \caption{A matrix $M$ with its $\rmin c_j$ and $\rmax c_j$ values}
    \label{MatFig}
\end{figure}

It will be convenient to index the rows and columns of our matrices by $0,1,\ldots,m$ for some $m\ge -1$ where $m=-1$ corresponds to the empty matrix.  In particular, we let $Z_m$ be the matrix with these coordinates which has all entries zero.  
Figure~\ref{MatFig} shows a matrix $M$ with rows and columns indexed by $[0,4]$.
We also let $E_{i,j}$ be the matrix whose only nonzero entry is a one in position $(i,j)$.  Note that we do not specify the dimensions of $E_{i,j}$ but instead choose them to be consistent with the dimensions of the other matrices involved in a given computation.  Given any matrix $M$ we let $c_j$ be the column vector consisting of its $j$th column.  If $c_j$ is not the zero vector then we define $\rmin c_j$ and $\rmax c_j$ to be the smallest and largest row coordinates, respectively, of a nonzero entry in $c_j$.  
For the matrix in  Figure~\ref{MatFig} the $\rmin c_j$ and $\rmax c_j$ values are listed above each column $c_j$.

If $d\ge1$ then a  {\em $d$-matrix} is a $[0,m]\times[0,m]$  matrix $M$ having the following properties.
\ben
\item[(M1)]  $M$ is upper triangular with entries $0$ and $1$.
\item[(M2)]  Between any two ones in the same column there are at least $d-1$ zeros.
\item[(M3)]  There are no zero columns and for all $j\in[m]$ we have 
$$
\rmax c_j > \rmin c_{j-1} - d.
$$
\een
For $n\ge0$ let
$$
\dMtx_n = \{ M \mid \text{$M$ is a $d$-matrix having $n$ ones}\}.
$$
It is easy to check that the matrix in Figure~\ref{MatFig} is in $\dMtx_8$ for $d=2$.

\begin{figure}
    \centering
$$
\barr{l}
\barr{cccccc}
 \left[\barr{ccccc} 
0&0&0&0&0\\ 
0&0&0&0&0\\ 
0&0&0&0&0\\ 
0&0&0&0&0\\ 
0&0&0&0&0
\earr
\right]
&
&
 \left[\barr{ccccc} 
1&0&0&0&0\\ 
0&0&0&0&0\\ 
0&0&0&0&0\\ 
0&0&0&0&0\\ 
0&0&0&0&0
\earr
\right]
&
&
 \left[\barr{ccccc} 
1&0&0&0&0\\ 
0&1&0&0&0\\ 
0&0&0&0&0\\ 
0&0&0&0&0\\ 
0&0&0&0&0
\earr
\right]
\\
M_0  & + E_{0,0} & M_1 & + E_{1,1} & M_2 & + E_{0,2}
\earr
\\
\\
\barr{cccccc}
 \left[\barr{ccccc} 
1&0&1&0&0\\ 
0&1&0&0&0\\ 
0&0&0&0&0\\ 
0&0&0&0&0\\ 
0&0&0&0&0
\earr
\right]
&
&
 \left[\barr{ccccc} 
1&0&1&0&0\\ 
0&1&0&0&0\\ 
0&0&0&0&0\\ 
0&0&0&1&0\\ 
0&0&0&0&0
\earr
\right]
&
&
 \left[\barr{ccccc} 
1&0&1&1&0\\ 
0&1&0&0&0\\ 
0&0&0&0&0\\ 
0&0&0&1&0\\ 
0&0&0&0&0
\earr
\right]
\\
M_3  & + E_{3,3} & M_4 & + E_{0,3} & M_5 & + E_{4,4}
\earr
\\
\\
\barr{cccccc}
 \left[\barr{ccccc} 
1&0&1&1&0\\ 
0&1&0&0&0\\ 
0&0&0&0&0\\ 
0&0&0&1&0\\ 
0&0&0&0&1
\earr
\right]
&
&
 \left[\barr{ccccc} 
1&0&1&1&0\\ 
0&1&0&0&0\\ 
0&0&0&0&1\\ 
0&0&0&1&0\\ 
0&0&0&0&1
\earr
\right]
&
&
 \left[\barr{ccccc} 
1&0&1&1&1\\ 
0&1&0&0&0\\ 
0&0&0&0&1\\ 
0&0&0&1&0\\ 
0&0&0&0&1
\earr
\right]
\\
M_6  & + E_{4,2} & M_7 & + E_{4,0} & M_8=\mx(\al)
\earr
\earr
$$
    \caption{Computing $\mx(0, 1, 0, 3, 0, 4, 2, 0)$}
    \label{mxFig}
\end{figure}

In order to construct our bijection $\mx:\dA_n\ra\dMtx_n$ we will need a certain factorization of $d$-ascent sequences.  Given  
$\al=a_1 a_2\ldots a_n\in \dA_n$, 
its {\em $d$-ascent factorization} is the concatenation 
\begin{equation}
\label{deFac}
\al=\de_0\de_1\ldots\de_m ,
\end{equation}
 where the factors $\de_i$ are obtained by dividing $\al$ after each $d$-ascent.  Equivalently, the $\de_i$ are the maximal factors of $\al$ containing no $d$-ascent.  Note that we start our indexing of the factors with $0$.
For example, if $d=2$ then $\al= 0, 1, 0, 3, 0, 4, 2, 0$ has factorization
$$
\al = 0 \cdot 1 \cdot 0 \cdot 3, 0 \cdot 4, 2, 0.
$$
where the  centered dots separate the factors.  Equivalently
$$
\de_0 = 0,\ \de_1 = 1,\ \de_2 = 0,\ \de_3 = 3, 0, \text{ and } \de_4= 4, 2, 0.
$$

Now given $\al=a_1\ldots a_n\in\dA_n$ expressed as in~\eqref{deFac} we will apply the map $\mx$ as follows.
Note that $m$ is the last index in the factorization and $n=|\al|$.
Construct a sequence of matrices
$$
Z_m = M_0, M_1,\ldots, M_n = \mx(\al)
$$
where 
\beq
\label{M_k}
M_k = M_{k-1} + E_{a_k,j}
\eeq
if $a_k$ is in the factor $\de_j$. 
This construction applied to the example $\al$ from the previous paragraph is shown in Figure~\ref{mxFig}, where each matrix appears below its label.

\bth
\label{MtxThm}
For all $d\ge 1$ and $n\ge0$ the map $\mx:\dA_n\ra\dMtx_n$ is a bijection.  Consequently
$$
\#\dA_n = \#\dMtx_n.
$$
\eth
\bprf
We first prove that $\mx$ is well defined in that $M_n\in\dMtx_n$.  
To ensure triangularity, we need to have $a_k\le j$ in equation~\eqref{M_k}.  Since $a_k\in \de_j$ and $\de_j$ is decreasing, it suffices to prove this inequality when $a_k$ is the first element of $\de_j$.  But the assumption about $a_k$ and $\de_j$ implies that the length $k-1$ prefix $\al_{k-1} = \de_0\de_1\ldots\de_{j-1}$.  
Now using the definition of a $d$-ascent sequence yields
$$
a_k \le \dasc \al_{k-1} + 1 = (j-1)+1 = j.
$$
So the first half of condition (M1) holds

We now show simultaneously that (M2) and the second half of (M1) are true.
If the ones added for $a_k$ and $a_{k+1}$ are in the same column then we must have
both in the same $\de_j$ for some $j$.
Since $\de_j$ has no $d$-ascents, it must be that  $a_{k+1}\le a_k-d$.
But since the elements of $\al$ give the row indices for the added ones, there must be at least $d-1$ zeros between the  ones corresponding to the $k$th and $(k+1)$st positions which is (M2).  Since
$d\ge1$ by assumption, we will never add two ones to the same position which finishes the proof of (M1).

 For condition (M3), since each $\de_j$ has at least one element, there are no zero columns.  Also, we have $\rmin c_{j-1}=a_{k-1}$  and $\rmax c_j=a_k$ where $a_{k-1}$ and $a_k$ are as  in the first  paragraph.
But then, by the definition of the factorization, there is a $d$-ascent from $a_{k-1}$ to $a_k$ and so
$$
\rmax c_j = a_k > a_{k-1}-d =\rmin c_{j-1} - d.
$$

To complete the proof of being well defined, we need to know that $\mx(\al)$ has exactly $n$ ones.  But an $E_{i,j}$ was added $n$ times and we already showed that $M_n$ is a $0$-$1$ matrix.  So $n$ must be the number of ones.

It is easy to construct a step-by-step inverse to $\mx$.  It follows that  this map is bijective and so we are done.
\eprf


\section{Matchings}
\label{matc}

We will construct our map between $\dA_n$ and certain matchings by restricting a bijection of Claesson and Linusson~\cite{CL:n!m} on inversion sequences.

Call a sequence of nonnegative integers $\al=a_1a_2\ldots a_n$
an {\em inversion sequence} if 
\begin{equation}
\label{inv}
 0\le a_k < k   
\end{equation}
for all $k\in[n]$.  The name derives from the fact that such sequences are in natural bijection with permutations in the symmetric group $\fS_n$ using the inversion statistic.
Let
$$
\I_n = \{\al=
a_1 a_2 \ldots a_n
\mid
\text{$\al$ is an inversion sequence}\}.
$$
For example,
$$
\I_3=\{000,001,002,010,011,012\}.
$$
Note that any $d$-ascent sequence is an inversion sequence:
Setting $k=1$ in~\eqref{inv} immediately gives $a_1=0$.
And for $k\ge2$, the maximum number of $d$-ascents in $\al_{k-1}$ is $k-2$ so that 
$$
a_k \le 1 +\dasc \al_{k-1} \le k-1.
$$
It follows that $\dA_n\sbe\I_n$ for all $d,n\ge0$.

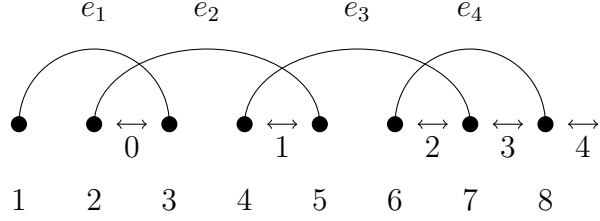
\begin{figure}
    \centering
\begin{tikzpicture}
\foreach \i in {1, 2, ..., 8} {
\filldraw(\i,1) circle(.1);
\draw(\i,0) node{$\i$};
};
\draw(3,1) arc (0:180:1);
\draw(5,1) arc (0:180:1.5 and 1);
\draw(7,1) arc (0:180:1.5 and 1);
\draw(8,1) arc (0:180:1);
\draw(2,2.5) node{$e_1$};
\draw(3.5,2.5) node{$e_2$};
\draw(5.5,2.5) node{$e_3$};
\draw(7,2.5) node{$e_4$};
\foreach \i in {2,4,6,7,8} {
\draw[<->] (\i+.3,1)--(\i+.7,1);
};
\draw(2.5,.7) node{$0$};
\draw(4.5,.7) node{$1$};
\draw(6.5,.7) node{$2$};
\draw(7.5,.7) node{$3$};
\draw(8.5,.7) node{$4$};
\end{tikzpicture}
    \caption{A matching and its active spaces}
    \label{mhFig}
\end{figure}

A {\em (perfect) matching}, $m$, is a graph on the vertex set $[2n]$ with a set of $n$ edges $e=ij$ no two  of which share a vertex.  We will always write our edges $ij$ so that $i<j$ and label them 
$e_1 = i_1 j_1, e_2 =i_2 j_2,\ldots, e_n = i_n j_n$ so that
$$
1< j_1< j_2 < \ldots < j_n = 2n.
$$
Figure~\ref{mhFig} displays a matching $m$ on $[8]$ with edges 
$$
e_1 = 13,\ e_2 = 25,\ e_3 = 47,\ \text{ and } e_4 = 68.
$$

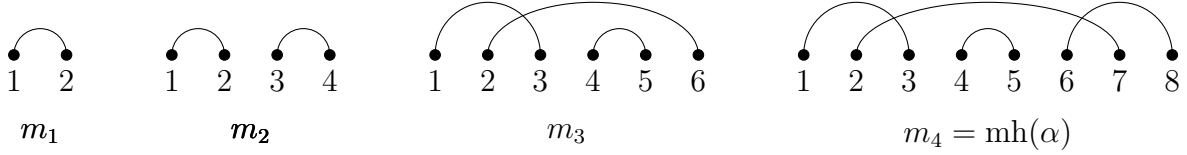
\begin{figure}
    \centering
\begin{tikzpicture}[scale=.7]
\foreach \i in {1, 2} {
\filldraw(\i,1) circle(.1);
\draw(\i,.5) node{$\i$};
\draw(1.5,-.5) node{$m_1$};
};
\draw(2,1) arc(0:180:.5);
\begin{scope}[shift={(3,0)}]
\foreach \i in {1, 2, ...,4} {
\filldraw(\i,1) circle(.1);
\draw(\i,.5) node{$\i$};
\draw(2.5,-.5) node{$m_2$};
};
\draw(2,1) arc(0:180:.5);
\draw(4,1) arc(0:180:.5);
\end{scope}
\begin{scope}[shift={(8,0)}]
\foreach \i in {1, 2, ...,6} {
\filldraw(\i,1) circle(.1);
\draw(\i,.5) node{$\i$};
};
\draw(3,1) arc(0:180:1);
\draw(5,1) arc(0:180:.5);
\draw(6,1)  arc(0:180:2 and 1);
\draw(3.5,-.5) node{$m_3$};
\end{scope}
\begin{scope}[shift={(15,0)}]
\foreach \i in {1, 2, ...,8} {
\filldraw(\i,1) circle(.1);
\draw(\i,.5) node{$\i$};
};
\draw(3,1) arc(0:180:1);
\draw(5,1) arc(0:180:.5);
\draw(7,1)  arc(0:180:2.5 and 1);
\draw(8,1) arc(0:180:1);
\draw(4.5,-.5) node{$m_4=\mh(\al)$};
\end{scope}
\end{tikzpicture}
    \caption{Computing $\mh(0,1,0,2)$}
    \label{mhSeq}
\end{figure}

The Claesson-Linusson bijection, $\mh$, is built on the notion of an active space.  The {\em active spaces} of $m$ are the spaces just before $j_1, j_2,\ldots, j_n$ together with the space after $2n$.
They will be labeled $0,1,\ldots,n$ from left to right.
In Figure~\ref{mhFig}, the active spaces are indicated using double-headed arrows with their labels below.
Given $\al=a_1a_2\ldots a_n\in\I_n$ we construct a sequence of matchings
$$
\emp = m_0, m_1,\ldots, m_n=\mh(\al)
$$
where $m_k$ is a matching on $[2k]$ for $0\le k\le n$.   To obtain $m_k$ from $m_{k-1}$, insert new vertices in active spaces $a_k$ and $k-1$ and connect them by an edge $e_k$
as well as renumbering the vertices $1,2,\ldots,2k$ from left to right.  For example, if $\al=0102$ then the sequence of nonempty matchings is given in Figure~\ref{mhSeq}.

To describe the image of the $\mh$ map, define a 
{\em nesting} in a matching $m$ to be a pair of edges of the form $ij$ and $kl$ where $i<k<l<j$.  As an example,
in Figure~\ref{mhSeq} the edges $26$ and $45$ form a nesting in $m_3$.  A {\em left nesting} is a nesting with $k=i+1$, and a {\em right nesting} is one with $j=l+1$.
The nesting in the previous example is a right nesting, but not a left nesting.
Let 
$$
\M_n =\{m \mid \text{$m$ is a matching on $[2n]$ with no left nestings}\}.
$$
Also, given an edge $e=ij$ in a matching we say that edge $kl$ is {\em to its left} if $l<i$.  We let
$$
\lef(e) = \text{ number of edges to the left of $e$}.
$$
We include the following result and proof of Claesson and Linsson to motivate what follows for $d$-ascent sequences.
\bth[\cite{CL:n!m}]
The map $\mh$ is a bijection $\mh:\I_n\ra\M_n$.
\eth
\bprf
The function is well defined since if $i_k j_k$ is the edge added to form $m_k$ then, by definition of the active spaces, $i_k + 1$ is the right endpoint of an edge.  So no left nesting can be formed.

To see that this is a bijection, one can construct the inverse by sending a matching $m$ with vertices $[2n]$ to the inversion table $\al=a_1 a_2\ldots a_n$ where
\beq
\label{lef}
a_k =\lef(e_k)
\eeq
for $k\in [n]$.
\eprf

Since every $d$-descent sequence is an inversion table, we can get a bijection with matchings by restricting $\mh$ to $\dA_n$ and characterizing the image.
In view of equation~\eqref{lef} and the definition of $\dAsc\al$ we say that $k$ is a {\em $d$-displacement} of a matching $m$ on $[2n]$ if
$$
\lef(e_{k+1}) > \lef(e_k)-d,
$$
and let 
$$
\ddis m=\#\{k\in[n-1] \mid \text{$k$ is a $d$-displacement of $m$}\}.
$$
We can now finally define the matchings which are in bijection with $d$-ascent sequences.  The set of 
{\em difference $d$ matchings} is
$$
\dMch_n = \{ m\in\M_n \mid \text{ $\lef(m_{k+1}) \le 1+ \ddis(m_k)$ for $k\in[n-1]$}\}
$$
where $m_k$ is the restriction of $m$ to the arcs $e_1,\ldots,e_k$ gotten by removing all vertices which are not endpoints of those edges and relabeling the remaining vertices.
The next result follows easily from the previous theorem, so its proof is omitted.
\bth
For all  $d,n\ge0$, the map $\mh$ restricts to a bijection $\mh:\dA_n\ra \dMch_n$.\hqed
\eth


\section{Permutations}
\label{per}

Bousquet-M\'elou et al~\cite{BMCDK:2fp} initiated the study of bivincular patterns and showed that ascent sequences are in bijection with permutations avoiding a bivincular pattern of length 3.  A bijection between weak ascent sequences and permutations avoiding a bivincular pattern of length 4 was found in~\cite{BCD:was}.
In this section, we will show that there is always an injection from $d$-ascent sequences into a corresponding set of permutations avoiding a pattern of length $d+3$, although it is not surjective in general.  In all cases, the technique used is based on the method of generating trees and, in particular, of active sites in permutations.

Recall that we are using $\fS_n$ for the symmetric group of all permutations $\pi=p_1 p_2\ldots p_n$ of $[n]$ written in one-line notation.  
Note, again, the use of Greek letters for permutations and Roman ones for their elements.
We say that $\pi$ contains the {\em classical pattern} $\si= s_1 s_2\ldots s_k\in\fS_k$ if there is a subsequence 
$\kappa=p_{i_1} p_{i_2} \ldots p_{i_k}$ of $\pi$ whose elements are in the same relative order as those of $\si$.
In this case $\ka$ is called a {\em copy} of $\si$.
For example, a copy of the pattern $231$ would be a subsequence $\ka=bca$ with $a<b<c$.
So $\pi=643512$ contains four copies of $231$ namely
$451$, $452$, $351$, and $352$.  To contain a {\em bivincular pattern} $\si$, certain pairs of elements of the copy $\ka$ must be adjacent in $\pi$ and others must be adjacent as integers.  In the first case we put a vertical bar between the elements of $\si$, and in the second we put a bar over the smaller of the two integers.
To illustrate, a copy $\ka =bca$ of the bivincular pattern
$\si'=2|3\ol{1}$ in a permutation $\pi$ would have $a<b<c$ with $b,c$ adjacent in $\pi$ and $b=a+1$.
Using $\pi$ as in the previous example, only one of the subsequences listed above is a copy of $\si'$, namely $352$.  Note that classical patterns are a special case of bivincular ones where there are no bars or overlines, so henceforth we will just use the term pattern. 
We say that $\pi$ {\em avoids} the pattern $\si$ if it does not contain $\si$ and let
$$
\Av_n(\si)=\{\pi\in\fS_n \mid \text{$\pi$ avoids $\si$}\}.
$$
\bth[\cite{BMCDK:2fp,BCD:was}]
Let
$$
\si_3 = 2|3\ol{1} \text{ and } \si_4 = 3|41\ol{2}.
$$
Then for all $n\ge0$ there are bijections $\A_n\leftrightarrow\Av_n(\si_3)$ and $\wA_n\leftrightarrow\Av_n(\si_4)$.  Thus
$$
\#\A_n=\#\Av_n(\si_3) \text{ and } \#\wA_n=\#\Av_n(\si_4).
$$
\eth

\begin{table}
    \centering
    \begin{tabular}{c||cccccccc}
       $n$          & 1 & 2 & 3 & 4 & 5    & 6      & 7    & 8      \\
       \hline\hline
     $\#\dA_n$      & 1 & 2 & 6 & 24& 118  & 682    & 4506 & 33376  \\
     \hline
     $\#\Av(\si_5)$ &1  & 2 & 6 & 24 & 119  & 699   & 4721 & 35904 
    \end{tabular}
    \caption{Comparison of $\#\dA_n$ when $d=2$ and $\#\Av(\si_5)$}
    \label{SiTab}
\end{table}

Consider the bivincular pattern
$$
\si_d = (d-1) | d 1 2 \ldots \ol{ (d-2)}.
$$
The main theorem of this section is as follows.
\bth
\label{AvSi}
For all $d,n\ge 0$ we have
$$
\#\dA_n \le \#\Av_n(\si_{d+3}).
$$
\eth
We note that the inequality in this theorem can be strict.  For example, when $d=2$ the cardinalities of
of $\dA_n$ and $\Av_n(\si_5)$ are shown in Table~\ref{SiTab} and are not always equal.

To prove Theorem~\ref{AvSi}, we will need the notion of an active site in a permutation.  Fix a pattern $\si$ and suppose that 
$\pi=p_1 p_2\ldots p_n\in\Av_n(\si)$.  Consider the $n+1$ spaces of $\pi$ consisting of the space before $\pi_1$, the space after $\pi_n$, and the $n-1$ spaces between adjacent elements of $\pi$.  Call a space an {\em active space} 
or {\em active site}
if insertion of $n+1$ in that space results in a permutation $\pi'$ which still avoids $\si$.  The rest of the sites/spaces are called
{\em inactive}.  We write
$$
\act\pi = \text{ number of active sites of $\pi$}
$$
and number the active sites of $\pi$ from left to right $0,1,\ldots,(\act\pi)-1$.  
Note that $\act\pi$ depends on the pattern $\si$ being avoided but this will always be clear from context.
The empty permutation $\emp$ has one active site labeled $0$.
For example, if $\pi=631245$ and $\si=\si_4=3|41\ol{2}$  then the numbering of its active sites is
$$
\rule{0pt}{5pt}_{\genfrac{}{}{0pt}{2}{\uparrow}{0}}6\ 3\ 1_{\genfrac{}{}{0pt}{2}{\uparrow}{1}}2_{\genfrac{}{}{0pt}{2}{\uparrow}{2}}4_{\genfrac{}{}{0pt}{2}{\uparrow}{3}}5_{\genfrac{}{}{0pt}{2}{\uparrow}{4}}.
$$
Note that the space after $6$ is not active because inserting $7$ there gives
the copy $\ka=6745$ of $\si$, and the space after $3$ is not active because of the creation of $\ka=3712$.

\begin{figure}
    \centering
$$
\barr{ccccccc}
\emp&
\rule{0pt}{5pt}_{\genfrac{}{}{0pt}{2}{\uparrow}{0}}1_{\genfrac{}{}{0pt}{2}{\uparrow}{1}}&
\rule{0pt}{5pt}_{\genfrac{}{}{0pt}{2}{\uparrow}{0}}1_{\genfrac{}{}{0pt}{2}{\uparrow}{1}}2_{\genfrac{}{}{0pt}{2}{\uparrow}{2}}&
\rule{0pt}{5pt}_{\genfrac{}{}{0pt}{2}{\uparrow}{0}}31_{\genfrac{}{}{0pt}{2}{\uparrow}{1}}2_{\genfrac{}{}{0pt}{2}{\uparrow}{2}}&
\rule{0pt}{5pt}_{\genfrac{}{}{0pt}{2}{\uparrow}{0}}31_{\genfrac{}{}{0pt}{2}{\uparrow}{1}}2_{\genfrac{}{}{0pt}{2}{\uparrow}{2}}4_{\genfrac{}{}{0pt}{2}{\uparrow}{3}}&
\rule{0pt}{5pt}_{\genfrac{}{}{0pt}{2}{\uparrow}{0}}31_{\genfrac{}{}{0pt}{2}{\uparrow}{1}}2_{\genfrac{}{}{0pt}{2}{\uparrow}{2}}5_{\genfrac{}{}{0pt}{2}{\uparrow}{3}}4_{\genfrac{}{}{0pt}{2}{\uparrow}{4}}&
631254
\\
\\
\pi_0 & \pi_1 & \pi_2 & \pi_3 & \pi_4 & \pi_5 & \pi_6 =\pe(\al)
\earr
$$
    \caption{Computing $\pe(010220)$ when $d=1$}
    \label{peFig}
\end{figure}

We now define an injection $\pe:\dA_n\ra\Av_n(\si_{d+3})$ as follows.  
Given $\al =a_1a_2\ldots a_n \in\dA_n$ we construct a sequence of 
permutations
$$
\emp =\pi_0,\pi_1,\ldots,\pi_n = \pe(\al)
$$
where $\pi_k$ is obtained from $\pi_{k-1}$ by inserting $k$ in the active space labeled $a_k$ for $1\le k\le n$.
As an example, suppose that $d=1$ so we are avoiding
$\si_{d+3}=\si_4 = 3|41\ol{2}$.
For the weak ascent sequence $\al=010220$ we would get the  sequence of permutation in Figure~\ref{peFig} where the active sites of each permutation have been labeled for convenience.

We must prove that this map is well defined in that there is a space of $\pi_{k-1}$ labeled by $a_k$.  This will be done using the next  lemma.
\begin{lem}
\label{LowUp}
For all $k\ge 1$ We have the following inequalities.
\ben
\item[(a)] If $\al=a_1 a_2\ldots a_{k-1}$ and $\be = \be a_k$ (concatenation) then
$$
\dasc\al\le \dasc\be \le \dasc\al + 1.
$$
\item[(b)] If $\si\in \Av_k(\si_{d+3})$ is obtained from $\pi\in \Av_{k-1}(\si_{d+3})$ by inserting $k$ in an active site then 
$$
\act\pi \le \act \si \le \act\pi + 1.
$$
\een
\end{lem}
\bprf
(a)  The  $d$-ascents of $\be$ are those in $\al$ plus possibly a new $d$-ascent  if $a_k>a_{k-1} - d$.  The result follows.

\smallskip

(b)  Call the sites  of $\pi$  which remain adjacent to the same elements in $\si$ {\em common}.  We first claim that common sites remain active or inactive in passing from $\pi$ to 
$\si$.  If a common site is inactive in $\pi$ then insertion of $k$ there forms a copy $\ka$ of $\si_{d+3}$.  But then insertion of $k+1$ in the same site forms a copy of $\si_{d+3}$ obtained by replacing $k$ by $k+1$ in $\ka$.  Now assume, towards a contradiction, that an active common site in $\pi$ becomes inactive in $\si$.  So there is a copy $\ka'$ of $\si_{d+3}$ when $k+1$ is inserted in the site of $\si$. That copy must contain $k+1$ acting as the $d+3$.  Also $k$ is not in $\ka'$, since if it were then it would be acting as the $d+2$.  But that implies that the site where $k+1$ was inserted is adjacent to $k$ and so not  common.  
Thus  inserting $k$ into $\pi$ in this common site produces a copy of $\si_{d+3}$ obtained by replacing the $k+1$ in $\ka'$ with $k$, which contradicts the fact that the site is active in $\pi$.

Now consider the new site of $\si$ just before and after $k$.  One can show in a manner similar to the previous paragraph that the site before $k$ will always be active.
So $\act\si=\act\pi+1$ or $\act\si=\act\pi$ depending on whether the site after $k$ is active or not, respectively.  This is what we wished to prove.
\eprf

We can now show that insertion  is well defined.
\bpr
\label{WellDef}
During the construction of $\pe(a_1 a_2\ldots a_n)$, there is  a site in $\pi_{k-1}$ labeled $a_k$ for $1\le k\le n$.
\epr
\bprf
By the way that sites are labeled, we must show that
$$
a_k \le \act \pi_{k-1} -1
$$
for all $k\in[n]$.  This is equivalent to showing that
\beq
\label{DascAct}
\dasc\al_k \le \act\pi_k- 2
\eeq
for all $k\in[n]$.  To see this, note that since $\al$ is a $d$-ascent sequence, the last inequlity implies
$$
a_k \le \dasc \al_{k-1} + 1 \le (\act\pi_{k-1} -2) + 1 =\act\pi_{k-1}-1
$$
with the converse being similar.

We will prove the inequality~\eqref{DascAct} by induction on $k$ where the  case $k=1$  is trivial.  Assume the inequality for $k$.  There are three cases depending on the value of $a_{k+1}$.

For the first case, suppose $a_{k+1}\le a_k-d$.  So $\dasc\al_{k+1}=\dasc\al_k$.  But by the previous lemma, the right side of~\eqref{DascAct} stays the same or goes up by one when passing from $\pi_k$ to $\pi_{k+1}$.  So it continues to hold for $k+1$.

For the second case, suppose $a_{k+1}> a_k$.  Now $\dasc\al_{k+1} = \dasc\al_k + 1$.  So it suffices to show that $\act\pi_{k+1} = \act\pi_k + 1$.  By the proof of Lemma~\ref{LowUp}, this will follow if we can show that the site after $k+1$ is active in $\pi_{k+1}$.  Since $a_{k+1}> a_k$, the site where $k+1$ is inserted is to the right of $k$.
Suppose, towards a contradiction, that inserting $k+2$ directly to the right of $k+1$ forms a copy of $\si_{d+3}$.  Then $k+1$ and $k+2$ take the roles of $d+2$ and $d+3$, respectively.  This forces $k$ to be the $d+1$ and so to be the right of $k+1$.  This contradiction shows that the site after $k+1$ is active as desired.

The last case is when $a_k-d<a_{k+1}\le a_k$.  Again $\dasc\al_{k+1} = \dasc\al_k + 1$.  So, as in the previous paragraph, we have to show that the site after $k+1$ is active in $\pi_{k+1}$.
As before, we assume it is not, so  inserting $k+2$ would form a copy $ k+1,k+2, k_1, k_2, \ldots, k_d, k$ of $\si_{d+3}$ where $k_1<k_2<\ldots<k_d$.  
We claim that there is some $i\in[d-1]$ such that all the spaces between $k_i$ and $k_{i+1}$ are inactive.  For suppose this was not the case.  Now the active space before $k$ is labeled $a_k$ and there are at least $d-1$ active spaces between $k_1$ and $k_d$.  So $k+1$ was inserted in a space labeled at most $a_k-d$ which contradicts $a_{k+1}>a_k-d$.

Suppose $i$ is chosen so there are no active spaces between $k_i$ and $k_{i+1}$.  Since $k_i<k_{i+1}$ there must be two adjacent elements $p_j<p_{j+1}$ of $\pi_k$ among the elements in the
factor of $\pi_k$ from $k_i$ to $k_{i+1}$.  Since the space between $p_j$ and $p_{j+1}$ is not active, placing $k+1$ there results in a copy $\ka$ of $\si_{d+3}$ where $p_j$ and $k+1$ plays the roles of $d+2$ and $d+3$, respectively.  Since $p_j<p_{j+1}$, replacing $k+1$ in $\ka$ by $p_{j+1}$ gives a copy of $\si_{d+3}$ in $\pi_k$ where $p_j$ and $p_{j+1}$ play the roles of $d+2$ and $d+3$, respectively.
This contradicts the fact that $\pi_k$ avoids $\si_{d+3}$ and completes the proof.
\eprf

We can now demonstrate Theorem~\ref{AvSi} which will be an immediate consequence of the following result.
\bth
For all $d,n\ge0$ the map $\pe:\dA_n\ra \Av_n(\si_{d+3})$ is an injection.
\eth
\bprf
From the previous proposition we know that the insertion process will produce a sequence $\pi_0,\pi_1,\ldots,\pi_n=\pe(\al)$.  Furthermore, each $\pi_i$ avoids $\si_{d+3}$ because elements are inserted in active sites.  To show that this is an injection, we induct on $n$.  Suppose $\al=a_1\ldots a_n$ and $\be = b_1\ldots b_n$ are distinct $2$-ascent sequences.
If $\al_{n-1}\neq\be_{n-1}$ then by induction $\pe(\al_{n-1})\neq \pe(\be_{n-1})$.  And insertion of $n$ into two different permutations can not make them equal so
$\pe(\al)\neq \pe(\be)$.  The other possibility is $\al_{n-1}=\be_{n-1}$ which forces $a_n \neq b_n$.  So $n$ is inserted in different positions of $\pe(\al_{n-1})=\pe(\be_{n-1})$ and we again have that $\pe(\al)$ and $\pe(\be)$ are distinct.
\eprf

\section{Posets}
\label{pos}

\begin{figure}
    \centering
\begin{tikzpicture}
\filldraw(1,0) circle(.1);
\filldraw(1,1) circle(.1);
\filldraw(1,2) circle(.1);
\filldraw(-1,2) circle(.1);
\filldraw(3,2) circle(.1);
\filldraw(2,3) circle(.1);
\draw (1,0)--(1,2)--(2,3)--(3,2)  (-1,2)--(1,1);
\draw(-2.5,1.5) node{$P=$};
\draw(1.5,0) node{$1$};
\draw(1.5,1) node{$2$};
\draw(1.5,2) node{$3$};
\draw(-1.5,2) node{$6$};
\draw(3.5,2) node{$4$};
\draw(2,3.5) node{$5$};
\end{tikzpicture}
\hs{20pt}
\begin{tikzpicture}
\filldraw(1,0) circle(.1);
\filldraw(1,1) circle(.1);
\filldraw(1,2) circle(.1);
\filldraw(-1,2) circle(.1);
\filldraw(3,2) circle(.1);
\filldraw(2,3) circle(.1);
\draw (1,0)--(1,2)--(2,3)--(3,2)  (-1,2)--(1,1);
\draw(-2.5,1.5) node{$Q=$};
\draw(1.5,0) node{$1$};
\draw(1.5,1) node{$2$};
\draw(1.5,2) node{$3$};
\draw(-1.5,2) node{$4$};
\draw(3.5,2) node{$5$};
\draw(2,3.5) node{$6$};
\end{tikzpicture}
    \caption{Two posets}
    \label{PosFig}
\end{figure}
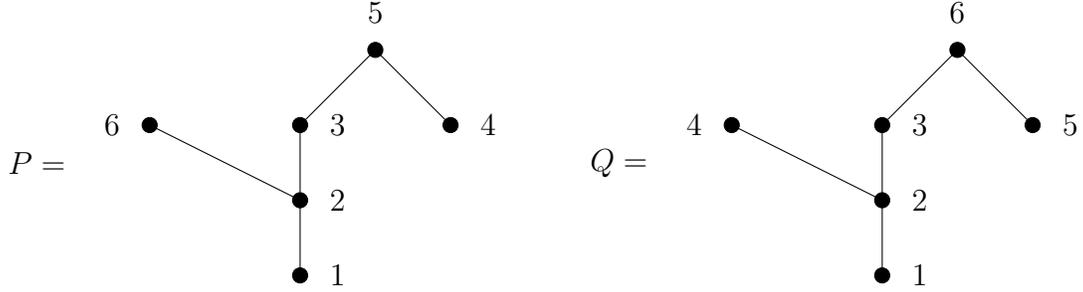

Let $(P,\le_P)$ be a poset (partially ordered set).
We say that $P$ is  {\em $(a+b)$-free} if it does not contain an (induced) subposet isomorphic to the disjoint union of an $a$-element chain and a $b$-element chain.  The poset $P$ whose Hasse diagram is on the left in Figure~\ref{PosFig} is not $(3+1)$-free because of the chains $1<_P 2 <_P 3$ and $4$.  But it is $(4+1)$-free because the only chain with four elements is $1<_P 2 <_P 3 <_P 5$ and the remaining two elements are related to members of that chain.
Ascent sequences are in bijection with unlabeled $(2+2)$-free posets as shown in~\cite{BMCDK:2fp}.  Claesson and Linusson~\cite{CL:n!m} introduced a family of labeled posets called factorial posets which are $(2+2)$-free.
Then in~\cite{BCD:was} a  bijective map was constructed between weak ascent sequences and factorial posets which do not contain a specially labeled 
$3+1$ subposet.
In this section we will introduce a specially labeled 
 $(d-1)+1$ poset $P_d$
and show that for all $d\ge1$ there is always an injection from $d$-ascent sequences to factorial posets which are special $P_{d+3}$-free.

Let $P$ be a poset with elements $[n]$ for some $n$.  We will use $<_P$ to denote the partial order on $P$ and $<$ for the total order on the integers.
Call $P$ {\em factorial} if it satisfies the following rule reminiscent of the transitive law:
\beq
\label{fact}
i<j \text{ and } j <_P k \implies i<_P k
\eeq
for all $i,j,k\in[n]$.  The reason they are called factorial is because there are $n!$ such posets on $[n]$.  As an example, it is easy to check that the poset $P$ on the left in Figure~\ref{PosFig} is factorial.  But the poset $Q$ on the right is not because $4<5$ and $5<_P 6$ but $4\not<_P 6$.  It is not hard to see that a factorial poset is {\em naturally labeled} in that $i<_P j$ implies $i<j$.

We define a {\em special  poset $P_d$} to be a disjoint union $(d-1)+1$ together with a labeling of the form
$$
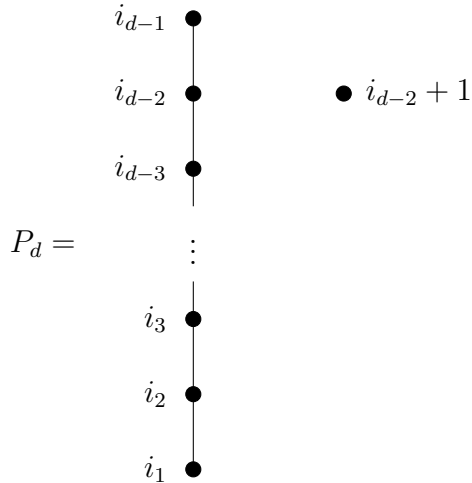
\begin{tikzpicture}
\filldraw(1,0) circle(.1);
\filldraw(1,1) circle(.1);
\filldraw(1,2) circle(.1);
\draw(1,3) node{$\vdots$};
\filldraw(1,4) circle(.1);
\filldraw(1,5) circle(.1);
\filldraw(1,6) circle(.1);
\filldraw(3,5) circle(.1);
\draw (1,0)--(1,2.5) (1,3.5)--(1,6);
\draw(-1,3) node{$P_d=$};
\draw(.5,0) node{$i_1$};
\draw(.5,1) node{$i_2$};
\draw(.5,2) node{$i_3$};
\draw(.3,4) node{$i_{d-3}$};
\draw(.3,5) node{$i_{d-2}$};
\draw(.3,6) node{$i_{d-1}$};
\draw(4,5) node{$i_{d-2}+1$};
\end{tikzpicture}
$$
where $i_1 < i_2 < \ldots <i_{d-1}$ and $i_{d-2}+1 < i_{d-1}$.
Call a labeled poset {\em special $P_d$-free} if it does not contain an induced subposet equal to $P_d$ for some choice of labels.  We will use the notation
$$
\Av_n(P_d) = \{P \mid \text{$P$ is a factorial poset on $[n]$ which is special $P_d$-free}\}.
$$
This should cause no confusion with the use of $\Av_n$ in the previous section since in that case it was applied to a permutation, not a poset.
To illustrate, note that even though the poset $P$ in Figure~\ref{PosFig} is not
$(3+1)$-free, it is special $P_4$-free.  This is because the unique copy of $3+1$ is
$1<_P 2 <_p 3$ and $4$, but the integer $4$ is not one more then the penultimate element of the $3$-element chain which is $2$.  This section's main result is as follows.
\bth
\label{AvP}
For all  $d\ge1$ and $n\ge0$ we have
$$
\#\dA_n \le \#\Av_n(P_{d+3}).
$$
\eth

\begin{table}
    \centering
    \begin{tabular}{c||cccccccc}
       $n$          & 1 & 2 & 3 & 4 & 5    & 6      & 7    & 8     \\
       \hline\hline
     $\#\dA_n$      & 1 & 2 & 6 & 24& 118  & 682    & 4506 & 33376  \\
     \hline
     $\#\Av(P_5)$ &1  & 2 & 6 & 24 & 119  & 700   & 4747 & 36370 
    \end{tabular}
    \caption{Comparison of $\#\dA_n$ when $d=2$ and $\#\Av(P_5)$}
    \label{PTab}
\end{table}

We note that for $d=1$ the above inequality is actually an equality as proved in~\cite{BCD:was}.  And, as in Theorem~\ref{AvSi}, the inequality can be strict for $d\ge2$ as  shown in Table~\ref{PTab}.

For $d\ge1$ we  define a map $\po:\dA_n\ra\Av_n(P_{d+3})$ as follows.
If $\al=a_1 a_2\ldots a_n$ then we construct a poset $P=\po(\al)$ on $[n]$ by letting, for each pair $i,j\in[n]$ of distinct elements,
\beq
\label{po(al)}
\text{$i<_P j$ if and only if $\dasc \al_i < a_j$.}
\eeq
For example, if $d=2$ and $\al=012032$ then $\po(\al)$ is the poset $P$ in Figure~\ref{PosFig}.  In particular, if $j=6$ then $a_6=2$ so that
$1<_P 6$ and $2<_P 6$
since $\dasc\al_1=0$ and $\dasc\al_2=1$ both of which are smaller than $a_6=2$.  However $6$ is not related to any $i>3$ since in those cases $\dasc\al_i\ge 2$.  We will first need a couple of lemmas to prove that $\po$ is well defined and an injection.

\begin{lem}
\label{aal}
Suppose $d\ge0$.  For any $d$-ascent sequence $\al=a_1 a_1\ldots a_n$ and any $k\in[n]$ we have
$$
a_k \le\dasc \al_k.
$$
\end{lem}
\bprf
Induct on $k$ where $k=1$ is trivial.  Assuming the result for $k-1$, we have two cases.
The first is if $a_k\le a_{k-1} -d$.  But then, by induction and the definition of a $d$-ascent,
$$
a_k \le a_{k-1} \le \dasc\al_{k-1} = \dasc\al_k.
$$
The other possibility is that $a_k>a_{k-1}-d$.  So, from the definition of a $d$-ascent sequence and induction, we have
$$
a_k \le \dasc\al_{k-1} + 1 = \dasc\al_k
$$
which completes the proof.
\eprf

We can now establish that the map $\po$ is well defined.
\begin{lem}
\label{PoWd}
Suppose $d\ge1$ and $\al\in\dA_n$.  Then $\po(\al)\in \Av_n(P_{d+3})$.
\end{lem}
\bprf
Let $\al=a_1 a_2\ldots a_n$ and $P=\po(\al)$.  
Since~\eqref{po(al)} only defines a relation on distinct elements, we can assume
reflexivity.  For antisymmetry, it suffices to show that if $i<_P j$ then
$i<j$ as then the law will follow from antisymmetry of the total order on integers.
Suppose, to the contrary, that $j<i$.  Using the previous lemma and the fact that $\dasc\al_k$ is a weakly increasing function of $k$ gives
$$
a_j \le \dasc \al_j\le \dasc \al_i .
$$
But this contradicts~\eqref{po(al)}.
For transitivity, suppose
$i<_P j$ and $j<_P k$.  By definition $\dasc \al_i< a_j$ and $\dasc \al_j < a_k$.  Now, using the previous sentence and Lemma~\ref{aal} again, we have
$$
\dasc \al_i < a_j \le \dasc \al_j < a_k.
$$
Thus $i<_P k$.

To check the factorial condition, assume $i<j$ and $j<_P k$.  This implies that $\dasc\al_i \le \dasc\al_j$ and $\dasc\al_j< a_k$.
By transitivity in the integers, $\dasc \al_i < a_k$ which is equivalent to $i<_P k$.

Finally, we need to show that $P$ is special $P_{d+3}$-free.
Take any chain
$$
C: i_1 <_P i_2 <_P \ldots < i_{d+1} <_P i_{d+2}
$$
in $P$ and, to simplify notation, let $i= i_1$, $j=i_{d+1}$, and $k = i_{d+2}$.
It suffices to show that $j+1$ must be related to some element of $C$.
There are now two cases.  First suppose that 
$a_{j+1}\le a_j -d$ so that $\dasc\al_{j+1}=\dasc\al_j$.
 Since $j<_p k$ we have $\dasc\al_j <a_k$.
Combining the inequalities gives $\dasc\al_{j+1}<a_k$ so that 
$j+1<_P k$,  thereby yielding the desired relation between $j+1$ and an element of $C$.

In the second case, suppose that $a_{j+1} > a_j - d$.
Note that since $d\ge1$, that element $i$ is different from $j$ and we will show $i<_P j+1$.
Now using the inequalities in $C$ and Lemma~\ref{aal} several more times
$$
\dasc \al_i =\dasc\al_{i_1} < a_{i_2}
\le \dasc\al_{i_2} < a_{i_3}\le \dasc\al_{i_3} < \ldots < a_{i_{d+1}} = a_j.
$$
It follows that $\dasc\al_i \le a_j-d$.
Combined with the inequality assumed for this case
we obtain $\dasc\al_i< a_{j+1}$.  This implies $i<_P j+1$ which completes the proof.
\eprf

We can now finish the proof that $\po$ is an injection.
\bth
For all $d\ge1$ and $n\ge0$ the map $\po:\dA_n\ra\Av(P_{d+3})$ is an injection.
\eth
\bprf
Since we have already proved that $\po$ is well defined, it just remains to prove that the map is injective.  It suffices to show that if $\al=a_1 a_2\ldots a_n$ with $\po(\al)=P$ then the sequence $\al$ is uniquely determined by $P$.  We will do this by induction on $n$ where, as usual, we skip the base case.  

By induction we can assume that $a_1,\ldots, a_{n-1}$ are determined.  There are now two cases.
First suppose that  $n$ is minimal in $P$.  We claim that this forces $a_n=0$.  For suppose $a_n>0$.
But then $a_n>0=\wasc\al_1$ and so $n>_P 1$, contradicting minimality  of $n$.

If $n$ is not minimal, then let
$$
m = \max\{i \mid i<_P n\}
$$
where the maximum is taken in the integers.
We claim that 
\beq
\label{analm}
a_n = \dasc\al_m + 1
\eeq
and so is uniquely determined.
Since $n>_P m$ we must have $a_n>\dasc\al_m$.
There are now two subcases depending on whether $m<n-1$ or $m=n-1$.
When $m<n-1$ we have $m+1<n$, and $m+1\not<_P n$ by definition of $m$.
So~\eqref{po(al)} implies that $\dasc\al_{m+1}\ge a_n$.  Thus
$$
1+\dasc\al_m \ge\dasc\al_{m+1}\ge a_n> \dasc\al_m
$$
which forces~\eqref{analm} to hold.  Now assume $m=n-1$ so that $a_n>\dasc\al_{n-1}$.  But the definition of a $d$-ascent sequence implies $a_n\le \dasc\al_{n-1}+1$.
This forces $a_n=\dasc\al_{n-1}+1$ as desired.
\eprf

\section{$d$-increasing $d$-ascent sequences}
\label{did}

Say that  a sequence of nonnegative integers $\al=a_1 a_2\ldots a_n$  is {\em $d$-increasing} if
$$
a_{k+1} > a_k - d
$$
for all $1\le k < n$.  Equivalently, 
\beq
\label{dinc}
\dAsc\al=[n-1].
\eeq
So a $0$-increasing sequence is strictly increasing while a $1$-increasing sequence is weakly increasing.
In this section we will study $d$-increasing $d$-ascent sequences.  In other words, we will consider the set
$$
\dI_n=\{\al\in\dA_n \mid \text{$\al$ is $d$-increasing}\}
$$
as well as the cardinalities
$$
\di_n =\#\dI_n.
$$
In particular, we will find a recursion for $\di_n$ and show that in the case $d=2$ these numbers count the rooted duplication trees introduced in~\cite{GHJM:ctd}.

It follows from equation~\eqref{dinc} that in a $d$-increasing sequence we have $\dasc\al_k = k-1$ for all $k$.
Combining this observation with the definitions of $d$-increasing and $d$-ascent sequence, we have that 
$\al=a_1 a_2\ldots a_n\in\dI_n$ if and only if
\ben
\item[(I1)] $a_1=0$, and
\item[(I2)] \label{I2} $a_k-d < a_{k+1} \le k$ for $k\in[n-1]$.
\een
Notice that if $d=0$ then there is a unique element in $\dI_n$, namely $\al=01\ldots n-1$.
When $d=1$, the  $\al\in\dI_n$ are weakly increasing sequences of nonnegative integers satisfying $a_k<k$ for all $k\in[n]$.  In this case $\di_n = C_n$, the $n$th Catalan number.  To see this, just interpret $a_k$ as the height  of the $k$th east step in a lattice path from the origin to $(n,n)$ using steps north and east and never going above the line $y=x$.  We will now derive a recursive formula for $\di_n$ for any $d\ge0$.  The sign-reversing involution proof which we use to obtain this result has wide applicability to sequence enumeration problems.  It was first used by Fr\"oberg~\cite{fro:dcp} and by Carlitz, Scoville and Vaughn~\cite{CSV:eps}.
For more information about sign-reversing involution proofs in general, see the text of Sagan~\cite{sag:aoc}.
\bth
\label{dinRec}
For $d\ge0$ and $n\ge1$ we have
$$
\di_n = \sum_{k\ge 1} (-1)^{k-1} \binom{n-kd+d}{k} \di_{n-k}.
$$
\eth
\bprf
When $d=0$ we have $\di_n=1$ for all $n$ and the identity is well known.  So we will assume that $d\ge1$.  
Putting everything on the left side of the equation, the desired formula can be rewritten
\beq
\label{dwSum}
\sum_{k\ge0} (-1)^k \binom{n-kd+d}{k} \di_{n-k} = 0.
\eeq
We will prove this by using a sign-reversing involution.  

Let $\cS$ be the set of all ordered pairs $(\al,S)$ consisting of  $\al=a_1 a_2\ldots a_{n-k}\in\dI_{n-k}$ for some $k\ge0$ and  a set
$S=\{s_1<s_2<\ldots<s_k\}\sbe\{0,1,\ldots,n-k\}$ such that 
\beq
\label{si}
s_{i+1}-s_i\ge d
\eeq
for all $1\le i<k$.
Notice that such $S$ are in bijection with unrestricted subsets 
$$
S'=\{s_1,s_2-d+1,s_3-2d+2,\ldots,s_k-(k-1)d+(k-1)\}
$$
of $\{0,1,\ldots,n-kd+k-1\}$ so that the number of possible $S$ is $\binom{n-kd+k}{k}$.
To each pair we associate a sign
$$
\sgn (\al,S) = (-1)^k.
$$
It follows directly from the definitions above that
$$
\sum_{(\al,S)\in\cS} \sgn (\al,S) = \sum_{k\ge0} (-1)^k \binom{n-kd+k}{k} \di_{n-k}.
$$
Thus it suffices to find a sign-reversing involution on $\cS$ with no fixed points.

Define $\io:\cS\ra\cS$ as follows.  Keeping the notation of the previous paragraph we let $\io(\al,S) = (\be,T)$ where there are two cases depending on $a_{n-k}$, the last element of $\al$, and $m$ which is the maximum element of $S$.   In particular, we define
$$
(\be,T) = \case{(\al m,\ S-\{m\})}{if the concatenation $\al m\in\dI_{n-k+1}$,}{(\al- a_{n-k},\ S\cup\{a_{n-k}\})}{otherwise,}
$$
where $\al-a_{n-k}$ is $\al$ with its last element removed.  We first need to check that this is well defined in that $(\be,T)\in\cS$.
First consider the case when  $\al m\in\dI_{n-k+1}$.  Then $\be=\al m$ is still a  $d$-increasing $d$-ascent sequence by definition and has length $n-k+1$.
Also $T=S-\{m\}$  is of size $k-1$.
By equation~\eqref{si} it
has largest element at most $n-k-1$ (in fact, at most $n-k-d$) and still satisfies that inequality for $0\le i<k-1$.
So  $(\be,T)\in\cS$ as desired.  Now consider the second case.  Removing an element from a  $d$-increasing $d$-ascent sequence does not change these characteristics, so 
$\be=\al- a_{n-k}\in\dI_{n-k-1}$.     By condition (I2) above we have $a_{n-k} \le n-k-1$.
So, $T\sbe\{0,1,\ldots,n-k-1\}$ which is the correct superset.  Finally, we claim that $T$ has cardinality $k+1$  and satisfies~\eqref{si} for $0\le i<k+1$. 
Because we are in the second case, it must be that $\al m$ is not a $d$-increasing $d$-ascent sequence.  Since $\al$ is  $d$-increasing, the number of $d$-ascents in $\al$ is $n-k-1$.
So concatenation of $\al$ with any number in $\{0,1,\ldots,n-k\}$ yields a $d$-ascent sequence.  Thus it must be that $\al m$ is not  $d$-increasing.  But this means that
$m\le a_{n-k}-d$.  Since $m$ is the maximum of $S$, it follows that $T = \{s_1 < s_2 <\ldots < m < a_{n-k}\}$ has $k+1$ elements  and that~\eqref{si} continues to hold.

We will be finished if we can show that $\io^2$ is the identity map since it clearly has no fixed points.
Given $(\al,S)\in\cS$ we let $\io(\al,S) = (\be,T)$.
Suppose first that $(\be,T)$ is obtained from $(\al,S)$ as in the first case of the definition of $\io$.
So $\be= \al m$ and $T=S-\{m\}$.  But now we can not remove the largest element of $T$ and append it to $\be$.  Indeed, $\be$ ends with $m$ and the largest element $l$ of $S-\{m\}$ satisfies $l\le m-d$.  So $\be l = \al m l$ is not $d$-increasing.  It follows that we must apply case two when computing $\io(\be,T)$ and will thus recover $(\al,S)$.
Now assume that $\io(\al,S)$ is computed using case two.  Then $\be=\al -a_{n-k}$ and, as shown in the previous paragraph, $a_{n-k}$ becomes the largest element of $T$.
Hence it is clearly possible to concatenate $\be$ and $a_{n-k}$ to reform $\al$ and the first case will be applied resulting in $\io(\be,T)=(\al,S)$,
\eprf

\bfig
\begin{tikzpicture}[scale=.8]
\draw(0,-1.7) node{};
\draw(1,0) grid (10,1);
\draw(0,.5) node{$\si'=$};
\draw(1.5,.5) node{$\si_1'$};
\draw(2.5,.5) node{$\si_2'$};
\draw(3.5,.5) node{$\si_3'$};
\draw(4.5,.5) node{$\si_4'$};
\draw(5.5,.5) node{$\si_5'$};
\draw(6.5,.5) node{$\si_6'$};
\draw(7.5,.5) node{$\si_7'$};
\draw(8.5,.5) node{$\si_8'$};
\draw(9.5,.5) node{$\si_9'$};
\draw(1,6) grid (7,7);
\draw(0,6.5) node{$\si=$};
\draw(1.5,6.5) node{$\si_1$};
\draw(2.5,6.5) node{$\si_2$};
\draw(3.5,6.5) node{$\si_3$};
\draw(4.5,6.5) node{$\si_4$};
\draw(5.5,6.5) node{$\si_5$};
\draw(6.5,6.5) node{$\si_6$};
\draw(0,3.5)  node{\rotatebox{270}{$\mapsto$}};
\draw(-.7,3.5)  node{$\phi_{2,3}$};
\draw[->] (1.5,5.5)--(1.5,1.5);
\draw[->] (2.5,5.5)--(2.5,1.5);
\draw[->] (3.5,5.5)--(3.5,1.5);
\draw[->] (4.5,5.5)--(4.5,1.5);
\draw[->] (2.5,5.5)--(5.5,1.5);
\draw[->] (3.5,5.5)--(6.5,1.5);
\draw[->] (4.5,5.5)--(7.5,1.5);
\draw[->] (5.5,5.5)--(8.5,1.5);
\draw[->] (6.5,5.5)--(9.5,1.5);
\draw[<->] (2.1,7.5)--(4.9,7.5);
\draw[<->] (5.1,7.5)--(6.9,7.5);
\draw(3.5,8) node {$r=3$};
\draw(6,8) node {$ a=2$};
\end{tikzpicture}
\hs{20pt}
\begin{tikzpicture}[scale=0.8]
\draw(-1,1) node{$T=$};
\node (n6) at (0,0) {6};
\node (nn1) at (1,1) {};
\node (n3) at (0,2) {3};
\node (nn2) at (3,1) {};
\node (nn3) at (4,0) {};
\node (nn4) at (4,2) {};
\node (n1) at (4,3) {1};
\node (n4) at (5,2) {4};
\node (n2) at (5,0) {2};
\node (n5) at (4,-1) {5};
\draw (nn1) -- (nn2) node[midway,above] {R};
\draw (n6) -- (nn1);
\draw (n3) -- (nn1);
\draw (nn2) -- (nn3);
\draw (nn3) -- (n2);
\draw (nn3) -- (n5);
\draw (nn2) -- (nn4);
\draw (nn4) -- (n1);
\draw (nn4) -- (n4);
\draw(2,-2)  node{\rotatebox{270}{$\mapsto$}};
\draw(2.7,-2)  node{$\phi_{2,3}$};
\begin{scope}[shift={(0,-6)}]
\draw(-1,1) node{$T'=$};
\node (n6) at (0,0) {9};
\node (nn1) at (1,1) {};
\node (n3) at (0,2) {};
\node (n3a) at (-0.25,3) {3};
\node (n3b) at (-1,2.25) {6};
\node (nn2) at (3,1) {};
\node (nn3) at (4,0) {};
\node (nn4) at (4,2) {};
\node (n1) at (4,3) {1};
\node (n4) at (5,2) {};
\node (n4a) at (6,2.5) {4};
\node (n4b) at (6,1.5) {7};
\node (n2) at (5,0) {};
\node (n2a) at (6,0.5) {2};
\node (n2b) at (6,-0.5) {5};
\node (n5) at (4,-1) {8};
\draw (nn1) -- (nn2) node[midway,above] {R};
\draw (n6) -- (nn1);
\draw (n3) -- (nn1);
\draw (n3) -- (n3a);
\draw (n3) -- (n3b);
\draw (nn2) -- (nn3);
\draw (nn3) -- (n2);
\draw (n2) -- (n2a);
\draw (n2) -- (n2b);
\draw (nn3) -- (n5);
\draw (nn2) -- (nn4);
\draw (nn4) -- (n1);
\draw (nn4) -- (n4);
\draw (n4)--(n4a);
\draw (n4)--(n4b);
\end{scope}
\end{tikzpicture}
\caption{A $(2,3)$-duplication}
 \label{phi23}
\efig

When $d=2$ there is a bijection between $2$-increasing $2$-ascent sequences and rooted duplication trees which we will now explain.  
Recall that we use the prefix t, rather than d, in our notation in the $d=2$ case.
A sequence of genes $\si=\si_1\si_2\ldots\si_n$ can evolve by having some factor $\si_i \si_{i+1}\ldots\si_{i+a}$ duplicate itself.
For example, on the left in Figure~\ref{phi23} the sequence $\si$ has the factor $\si_2\si_3\si_4$ duplicate itself as indicated by the arrows to form $\si'$.
The duplication process is specified by two parameters: $a\ge0$ which is the number of genes after the duplicated segment, and $r\ge1$ which is the number of genes in the duplicated segment.
This is called a {\em $(a,r)$-duplication} and we write $\si'=\phi_{a,r}(\si)$.
In the example of the figure, $a=2$, $r=3$ and $\si'=\phi_{2,3}(\si)$.

One can model a sequence of duplications using certain trees $T$.  Such $T$ have a rooted edge, $R$, indicating the beginning of the process.  They are also binary because a factor is doubled at each stage.
The leaves are bijectively labeled by $[n]$ for some $n$ and represent the genes $\si_1,\ldots,\si_n$.  The distance between two leaves in the tree weakly increases with each duplication.  If two leaves 
$\{\ell,m\}$ are siblings then we call them a {\em leaf pair}.  Suppose that tree $T$ corresponds to a sequence $\si=\si_1\ldots\si_n$ on which an $(a,r)$-duplication is performed to obtain $\si'$.  Then this is mirrored in a tree $T'=\phi_{a,r}(T)$ obtained from $T$ as follows.
\ben
\item  The leaves $\ell\in[n-a+1,n]$ are all renumbered as $\ell+r$.
\item  Each leaf $\ell\in[n-a-r+1,n-a]$ becomes the parent of a leaf pair $\{\ell,\ell+r\}$.
\een
If the sequence $\si$ in Figure~\ref{phi23} has resulted from a history whose duplication tree is $T$ as given at the  top-right, the resulting $T'=\phi_{2,3}(T)$  corresponding to $\si'$ is displayed on the bottom-right.

We can now formally define the set of rooted duplication trees, $\RDT$, inductively as follows.
\ben
\item  The tree
\bce
\begin{tikzpicture}[scale=0.6]
\node at (-2,0) {$T_2:=$};
\node (n1) at (0,0) {1};
\node (n2) at (3,0) {2};
\draw (n1) -- (n2) node[midway,above] {R};
\end{tikzpicture}
\ece
is in $\RDT$.
\item  If $T\in\RDT$ has $n$ leaves then $\phi_{a,r}(T)\in\RDT$ for any $a,r$ such that $a\ge 0$ and $1\le r\le n-a$.
\een
We also let 
$$
\RDT_n = \{T\in\RDT \mid \text{ $T$ has $n$ leaves}\}.
$$
Note that $T_2$ will always stand for the unique element in $\RDT_2$.

\bfig
\begin{center}
\footnotesize
\begin{tikzpicture}[scale=0.8]
\node (n6) at (0,0) {12};
\node (nn1) at (1,1) {};
\node (n3) at (0,2) {};
\node (n3a) at (-0.25,3) {5};
\node (n3b) at (-1,2.25) {8};
\node (nn2) at (3,1) {};
\node (nn3) at (4,0) {};
\node (nn4) at (4,2) {};
\node (nn5) at (4,3) {};
\node (nn6) at (5,3.5) {};
\node (n1) at (4,4) {1};
\node (new2) at (6,3.95) {2};
\node (new3) at (6,3.25) {3};
\node (n4) at (5,2) {};
\node (n4a) at (6,2.5) {6};
\node (n4b) at (6,1.5) {9};
\node (n2) at (5,0) {};
\node (n2a) at (6,0.5) {4};
\node (n2b) at (6,-0.5) {7};
\node (n5) at (4,-1) {};
\node (new10) at (3.5,-2) {10};
\node (new11) at (4.5,-2) {11};
\draw (nn1) -- (nn2) node[midway,above] {R};
\draw (n6) -- (nn1);
\draw (n3) -- (nn1);
\draw (n3) -- (n3a);
\draw (n3) -- (n3b);
\draw (nn2) -- (nn3);
\draw (nn3) -- (n2);
\draw (n2) -- (n2a);
\draw (n2) -- (n2b);
\draw (nn3) -- (n5);
\draw (n5) -- (new10);
\draw (n5) -- (new11);
\draw (nn2) -- (nn4);
\draw (nn4) -- (nn5);
\draw (nn5) -- (n1);
\draw (nn5) -- (nn6);
\draw (nn6) -- (new2);
\draw (nn6) -- (new3);
\draw (nn4) -- (n4);
\draw (n4)--(n4a);
\draw (n4)--(n4b);
\end{tikzpicture}
\end{center}
\caption{A tree in $\RDT_{12}$} 
\label{T12}
\efig

A complication with this definition of $\RDT$ is that the same tree can be constructed  by different applications of the  $\phi$-maps.  For example, $\phi_{2,1}(\phi_{0,1}(T_2))=\phi_{0,1}(\phi_{1,1}(T_2))$.
Gasceul et al.~\cite{GHJM:ctd} gave a canonical way to describe a given tree $T$ based on picking one of the most recent duplications which could lead to $T$.
Say that $T$ contains a {\em visible $(a,r)$-duplication event} if there is a tree $T'$ such that $\phi_{a,r}(T')=T$.  Note that the roles of $T$ and $T'$ have switched from the way they were used in  Figure~\ref{phi23}.  Having a visible $(a,r)$-duplication event for $T\in\RDT_n$ is equivalent to there being $a,r$ such that  $\{\ell,\ell-r\}$ are leaf pairs for 
$\ell\in[n-a-2r+1,n-a]$.  For example, suppose $T$ is as in Figure~\ref{T12}.  The visible duplication events in $T$ are as follows:
\begin{itemize}
\item $\{2,3\}$ is a leaf pair so $T$ contains a visible $(9,1)$-duplication event,
\item  $\{ 4,7\},\ \{5,8\},\ \{6,9\}$ are leaf pairs so  $T$ contains a visible $(3,3)$-duplication event,
\item $\{10,11\}$ is a leaf pair so $T$ contains a visible $(1,1)$-duplication event.
\end{itemize}
Importantly, if $T$ contains a visible $(a,r)$-duplication event then $a$ uniquely determines $r$ since leaf $n-a$ must have a sibling.
Moreover, if the label of that sibling is $m$ then $r=n-a-m$.
The {\em leftmost visible $(a,r)$-duplication event} for $T$ is the one with largest $a$.  By the observation just made, this is well defined.

\bfig
\begin{tikzpicture}[scale=0.7]
\draw(-1,1) node{$T_8=$};
\node (n6) at (0,0) {8};
\node (nn1) at (1,1) {};
\node (n3) at (0,2) {4};
\node (nn2) at (3,1) {};
\node (nn3) at (4,0) {};
\node (nn4) at (4,2) {};
\node (nn5) at (4,3) {1};
\node (n4) at (5,2) {5};
\node (n2) at (5,0) {};
\node (n2a) at (6,0.5) {2};
\node (n2b) at (6,-0.5) {3};
\node (n5) at (4,-1) {};
\node (new10) at (3.5,-2) {6};
\node (new11) at (4.5,-2) {7};
\draw (nn1) -- (nn2) node[midway,above] {R};
\draw (n6) -- (nn1);
\draw (n3) -- (nn1);
\draw (nn2) -- (nn3);
\draw (nn3) -- (n2);
\draw (n2) -- (n2a);
\draw (n2) -- (n2b);
\draw (nn3) -- (n5);
\draw (n5) -- (new10);
\draw (n5) -- (new11);
\draw (nn2) -- (nn4);
\draw (nn4) -- (nn5);
\draw (nn4) -- (n4);
\begin{scope}[shift={(11,0)}]
\draw(-1,1) node{$T_7=$};
\draw(-3,1) node{$\stackrel{\psi_{5,1}}{\mapsto}$};
\node (n6) at (0,0) {7};
\node (nn1) at (1,1) {};
\node (n3) at (0,2) {3};
\node (nn2) at (3,1) {};
\node (nn3) at (4,0) {};
\node (nn4) at (4,2) {};
\node (nn5) at (4,3) {1};
\node (n4) at (5,2) {4};
\node (n2) at (5,0) {2};
\node (n5) at (4,-1) {};
\node (new10) at (3.5,-2) {5};
\node (new11) at (4.5,-2) {6};
\draw (nn1) -- (nn2) node[midway,above] {R};
\draw (n6) -- (nn1);
\draw (n3) -- (nn1);
\draw (nn2) -- (nn3);
\draw (nn3) -- (n2);
\draw (nn3) -- (n5);
\draw (n5) -- (new10);
\draw (n5) -- (new11);
\draw (nn2) -- (nn4);
\draw (nn4) -- (nn5);
\draw (nn4) -- (n4);
\end{scope}
\end{tikzpicture}

\vs{10pt}

\begin{tikzpicture}[scale=0.7]
\draw(-1,1) node{$T_6=$};
\draw(-3,1) node{$\stackrel{\psi_{1,1}}{\mapsto}$};
\node (n6) at (0,0) {6};
\node (nn1) at (1,1) {};
\node (n3) at (0,2) {3};
\node (nn2) at (3,1) {};
\node (nn3) at (4,0) {};
\node (nn4) at (4,2) {};
\node (nn5) at (4,3) {1};
\node (n4) at (5,2) {4};
\node (n2) at (5,0) {2};
\node (n5) at (4,-1) {5};
\draw (nn1) -- (nn2) node[midway,above] {R};
\draw (n6) -- (nn1);
\draw (n3) -- (nn1);
\draw (nn2) -- (nn3);
\draw (nn3) -- (n2);
\draw (nn3) -- (n5);
\draw (nn2) -- (nn4);
\draw (nn4) -- (nn5);
\draw (nn4) -- (n4);
\begin{scope}[shift={(9,0)}]
\draw(0,1) node{$T_5=$};
\draw(-2,1) node{$\stackrel{\psi_{0,3}}{\mapsto}$};
\node (nn1) at (1,1) {5};
\node (nn2) at (3,1) {};
\node (nn3) at (4,0) {};
\node (nn4) at (4,2) {};
\node (nn5) at (4,3) {1};
\node (n4) at (5,2) {3};
\node (n2) at (5,0) {2};
\node (n5) at (4,-1) {4};
\draw (nn1) -- (nn2) node[midway,above] {R};
\draw (nn2) -- (nn3);
\draw (nn3) -- (n2);
\draw (nn3) -- (n5);
\draw (nn2) -- (nn4);
\draw (nn4) -- (nn5);
\draw (nn4) -- (n4);
\end{scope}
\end{tikzpicture}

\vs{10pt}

\begin{tikzpicture}[scale=0.7]
\draw(0,1) node{$T_4=$};
\draw(-2,1) node{$\stackrel{\psi_{1,2}}{\mapsto}$};
\node (nn1) at (1,1) {4};
\node (nn2) at (3,1) {};
\node (nn3) at (4,0) {3};
\node (nn4) at (4,2) {};
\node (nn5) at (4,3) {1};
\node (n4) at (5,2) {2};
\draw (nn1) -- (nn2) node[midway,above] {R};
\draw (nn2) -- (nn3);
\draw (nn2) -- (nn4);
\draw (nn4) -- (nn5);
\draw (nn4) -- (n4);
\begin{scope}[shift={(9,0)}]
\draw(0,1) node{$T_3=$};
\draw(-2,1) node{$\stackrel{\psi_{2,1}}{\mapsto}$};
\node (nn1) at (1,1) {3};
\node (nn2) at (3,1) {};
\node (nn3) at (4,0) {2};
\node (nn4) at (4,2) {1};
\draw (nn1) -- (nn2) node[midway,above] {R};
\draw (nn2) -- (nn3);
\draw (nn2) -- (nn4);
\end{scope}
\end{tikzpicture}

\vs{10pt}

\begin{tikzpicture}[scale=.7]
\draw(0,1) node{$T_2=$};
\draw(-2,1) node{$\stackrel{\psi_{1,1}}{\mapsto}$};
\draw(1,0) node{};
\node (nn1) at (1,1) {2};
\node (nn2) at (3,1) {1};
\draw (nn1) -- (nn2) node[midway,above] {R};
\end{tikzpicture}
\caption{A canonical reduction history}
\label{CrhFig}
\efig

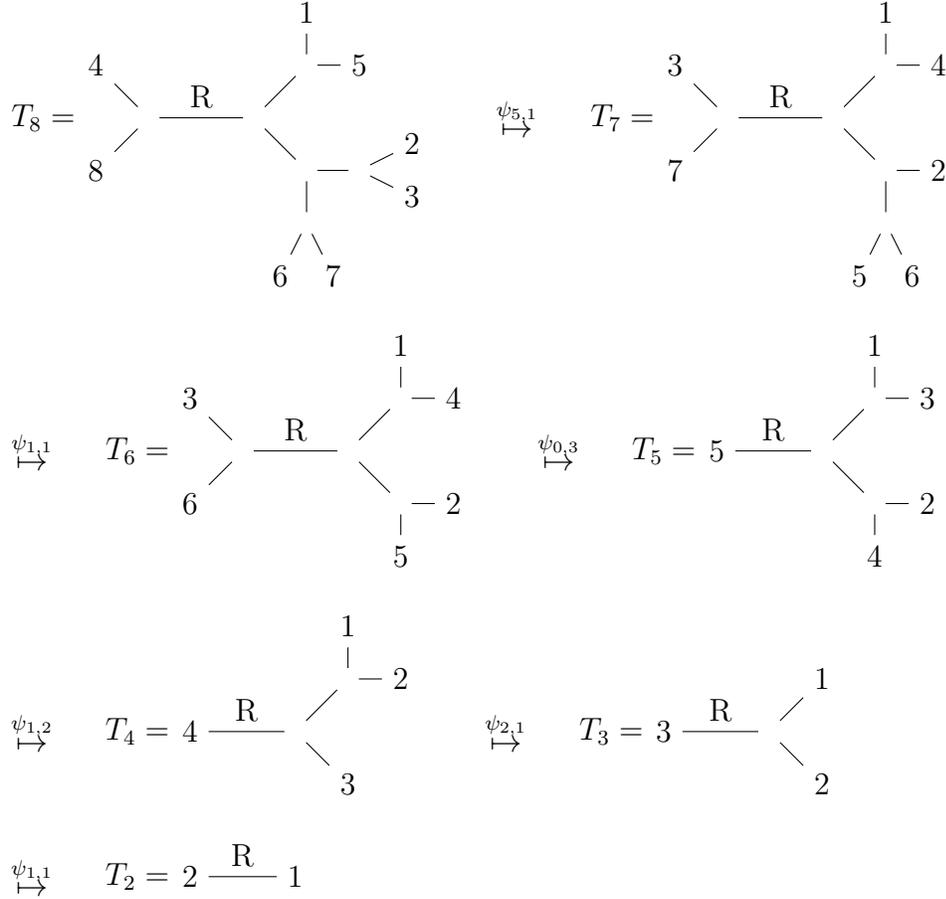
We now define the {\em canonical reduction history} of $T\in\RDT_n$ to be the sequence of trees
\beq
\label{crh}
T= T_n, T_{n-1}, T_{n-2},\ldots,T_2
\eeq
such that $T_k\in\RDT_k$ for $k\in[2,n]$ as follows.
If $T_k$, where $k\in[3,n]$, has leftmost duplication event with parameters $(a,r)=(a_k,r_k)$ then we will construct $T_{k-1}= \psi_{a,r}(T_k)$ by 
performing an operation which corresponds to
reducing the length of the duplicated string by one
and call $T_{k-1}$ the {\em canonical $(a,r)$-reduction} of $T_k$.
Formally, $T_{k-1}$ is obtained from $T_k$ by
\ben
\item[(1)]  The leaf pair $\{k-a,k-a-r\}$ is removed and the parent is labeled $k-a$.
\item[(2)]  Each leaf $\ell>k-a-r$ is relabeled $\ell-1$.
\een
So after step (1) 
the leaves are labeled by $[k]\setm\{k-a-r\}$ and then step (2) 
adjusts the labels to be in $[k-1]$.
We also define the {\em canonical reduction sequence} of $T$ to be
$$
\al(T)=(a_2,a_3,\ldots,a_n)
$$
where $a_2 {\bf:=} 0$.
Figure~\ref{CrhFig} shows the canonical reduction history for the tree in its upper-left corner.  So the tree $T=T_8$ has canonical reduction sequence
$$
\al(T) = (0, 1, 2, 1, 0, 1, 5).
$$

In order to prove that the map $T\mapsto\al(T)$ is a bijection from $\RDT_n$ to $\tI_{n-1}$ we will need the following lemma.
\begin{lem}
\label{aa'}
Suppose that in a canonical reduction history we have three consecutive trees $T$, $T'=\psi_{a,r}(T)$, and $T''=\psi_{a',r'}(T')$.
\ben
\item[(a)] If $r\ge2$ then 
$$
a'=a+1 \text{ and } r' = r-1.
$$
\item[(b)]  If $r=1$ then
$$
a' \le a.
$$
\een
\end{lem}
\bprf
(a)  Since $\psi_{a,r}$ was applied to $T$, that tree's  leftmost visible duplication  consisted of leaf pairs $\{\ell,\ell-r\}$ for 
$\ell\in[n-a-r+1,n-a]$ where $a$ is maximum.  In $T'$, the leaf pair $\{n-a,n-a-r\}$ was replaced by a leaf $n-a-1$ and the rest of the pairs
became 
$\{\ell-1,\ell-r\}$ for $\ell\in[n-a-r+1,n-a-1]$.  Since $r\ge2$ there is at least one such pair, and these pairs form a visible $(a+1,r-1)$-duplication event in $T'$.
So to finish the demonstration, we just need to show they are leftmost.

Suppose, towards a contradiction, that in $T'$ there was a visible $(a'',r'')$-duplication event  with $a''>a+1$.  But then $n-a''\le n-a-2$ and $n-a-2$ is in one of the pairs for the 
$(a,r)$-duplication event in $T$. So the $(a'',r'')$-duplication event in $T'$ must correspond to a duplication event in $T$ which is to the left of the $(a,r)$-duplication event.  This contradicts the fact
that the $(a,r)$-duplication event is leftmost in $T$.

(b)  The proof is similar to that of the second paragraph of the demonstration of (a).  So it is left to the reader.
\eprf

\bth
\label{RDTThm}
The map $T\mapsto\al(T)$ is a bijection $\RDT_n\ra\tI_{n-1}$.
\eth
\bprf
We first need to show that the map is well defined in that $\al(T)=(a_2,\ldots,a_n)\in\tI_{n-1}$.  We have $a_2=0$ by definition.
So we only need to verify that the inequalities in (I2) hold for $d=2$.  Taking into account that $\al(T)$ starts with index $2$, these inequalities translate
as 
$$
a_k - 1\le  a_{k+1}\le k-1.
$$
We prove this by induction on $k$ where the base case is easy.

Assume that the inequalities hold for $k-1$ and let the canonical reduction history trees of $T$ be as in~\eqref{crh}.
Also, adopt the notation of the previous lemma where $T=T_{k+1}$, $T'=T_k$, and $T''=T_{k-1}$.
So we wish to prove that 
$$
a'-1 \le a \le k-1.
$$
The largest $a$ can be is if $T$ has a single visible duplication event consisting of leaf pair $\{1,2\}$.  Since $T$ has $k+1$ leaves, this results in $a=(k+1)-2=k-1$ which is the desired upper bound.
Cases (a) and (b) of Lemma~\ref{aa'} immediately imply that $a\ge a'-1$ which  confirms the  lower bound.

To show we have a bijection we will construct the inverse.
Given $\al=(a_2,a_3,\ldots,a_n)\in\tI_{n-1}$ we will construct a sequence of trees $T_2,T_3,\ldots,T_n=T$ and nonnegative integers $r_2,r_3,\ldots,r_n$.
We start with $r_2=1$ and $T_2$ as the unique tree in $\RDT_2$.  Given $T_k$ and $r_k$ we construct $T_{k+1}$ and $r_{k+1}$ as follows.
First, we let
$$
r_{k+1} = \case{r_k+1}{if $a_{k+1} = a_k - 1$,}{1}{otherwise.}
$$
Letting $a=a_{k+1}$ and $r=r_{k+1}$, we create  $T_{k+1}$ from $T_k$ by
\ben
\item Each leaf $\ell>k-a-r$ is relabeled $\ell+1$.
\item The leaf labeled $k-a+1$ is given a sibling $k-a-r+1$.
\een
We leave the verification that this is well defined and the inverse of the map $T\mapsto\al(T)$ as a straightforward exercise.
\eprf

\section{Restricted growth functions}
\label{rgf}

The reader will have noticed that Theorem~\ref{AvP} is only true for $d\ge 1$ and so only applies starting with the poset $P_4$.  It is natural to ask if anything can be said about factorial posets which are special $P_3$-free, and it turns out that these are in natural bijection with restricted growth functions (RGFs).  These sequences are a way of modeling set partitions and have been widely studied~\cite{BKM:apr,CDDGGPS:rgf,kad:win,LF:1ar,MNV:lgc,mil:rgf1,mil:qrg,mil:rgf2,SV:mrg,wac:srg,WW:psn}.
We will also see how these sequences are related to pattern avoidance, matrices, and matchings.

A {\em restricted growth function (RGF)  of length $n$} is a sequence of nonnegative integers $\rho=r_1 r_2 \ldots r_n$ satisfying 
\ben
\item[(r1)] $r_1=0$, and 
\item[(r2)] $r_{k+1} \le 1+\max\rho_k$ for $k\in[n-1]$,
\een
where, as with $d$-ascent sequences, we let $\rho_k = r_1 r_2\ldots r_k$ and the maximum of a sequence is the maximum value of its elements as integers.
Since $\rho_{k+1} =\rho_k r_{k+1}$ and $\rho$ is an RGF,  it must be that for all $k<n$ we have
\beq
\label{k,k+1}
\max\rho_k \le \max\rho_{k+1}\le \max\rho_k + 1,
\eeq
where $\rho_{k+1}$ achieving its upper bound forces $r_{k+1}=\max\rho_{k+1}$.
Let
$$
\RGF_n = \{\rho \mid \text{$\rho$ is a restricted growth function of length $n$}\}.
$$
It is not hard to see that 
$$\RGF_n\sbe \A_n,
$$
the set of length $n$ ascent sequences.
Also, the RGFs of length $n$ are in natural bijection with the partitions of $[n]$ so that 
\beq
\label{RGFB}
\#\RGF_n = B_n,
\eeq
the $n$th Bell number.

We will show that there is a bijection $\rp:\RGF_n \ra\Av_n(P_3)$.
Given $\rho\in\RGF_n$ we define a poset $P=\rp(\rho)$ by
\beq
\label{rhor}
\text{$i<_P j$ if and only if $\max\rho_i < r_j$}.
\eeq
Note the similarity with the definition of $\po(\al)$~\eqref{po(al)}.
We must first show that this map is well defined.
\begin{lem}
If $\rho\in\RGF_n$ then $P=\rp(\rho)\in\Av_n(P_3)$.
\end{lem}
\bprf
The fact that $P$ is a poset is proved in much the same way as in Lemma~\ref{PoWd}, and so is omitted.
To show that $P$ is factorial, assume $i<j$ and $j<_P k$.  Since $\rho_i$ is a subsequence of $\rho_j$ we have $\max\rho_i\le\max\rho_j$.
 Also, $j<_P k$ implies $\max\rho_j<r_k$.  So, by transitivity, $\max\rho_i < r_k$ which is equivalent to $i<_P j$.

We show that $P$ avoids a special $P_3$ by induction on $n$, where the statement is obvious for $n\le 2$.  Assume, that for RGFs of length $n-1$ their image under $\rp$ is special $P_3$-free.  Now consider an RGF $\rho$ of length $n$ and assume, towards a contradiction, that $P=\rp(\rho)$ contains a special $P_3$ labeled as
$$
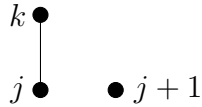
\begin{tikzpicture}
\filldraw(1,1) circle(.1);
\filldraw(1,2) circle(.1);
\filldraw(2,1) circle(.1);
\draw (1,1)--(1,2);
\draw(.7,1) node{$j$};
\draw(.7,2) node{$k$};
\draw(2.7,1) node{$j+1$};
\end{tikzpicture}
$$
where $j+1<k$.  
Now the prefix $\rho'=\rho_{n-1}$ of $\rho$ is of length $n-1$.  So, by induction, $P'=\rp(\rho')$
is special $P_3$-free.  But by~\eqref{rhor}, the poset $P'$ is the restriction of $P$ to the elements of the set $[n-1]$.  So the only way that $P'$ could avoid $P_3$ while $P$ contains a copy would be if that copy includes $n$.  Since $k$ is the largest integer in our copy, it must be that $k=n$.

Now $j<_P n$ and $j+1 \nless_P n$ implies $\max\rho_j < r_n\le\max\rho_{j+1}$.
   Combining the previous sentence with equation~\eqref{k,k+1} and the comment which follows it, we must  have  $r_{j+1}=\max\rho_{j+1} =\max\rho_j + 1$.
But  $j\nless_P j+1$  and so $\max\rho_j \ge r_{j+1}$ which contradicts the previous sentence.
\eprf

In order to provide an inverse for $\rp$, we will need the following lemmas.  For $k$ in a factorial poset $P$ we  will  use the notation
$$
M(k) =\max\{i \in P \mid i<_P k\},
$$
where the maximum is taken in the integers.  It follows easily from properties of factorial posets that if $M(k)=m$ then
\beq
\label{MkInt}
\{i \in P \mid i<_P k\} = [m].
\eeq
\begin{lem}
If $\rho=r_1 r_2\ldots r_n$ and $P=\rp(\rho)$ then for all $k\in[n]$ we have
\beq
\label{rp-1}
r_k=\case{0}{if $k$ is minimal in $P$,}{\max\rho_{M(k)} + 1}{otherwise.}
\eeq
\end{lem}
\bprf
Suppose first that $k$ is  minimal in $P$.  If $k=1$ then $r_1=0$ by definition of an RGF.  If $k>1$ then we must have $\max\rho_j\ge r_k$ for all $j<k$.  Letting $j=1$ this becomes $0=\max\rho_1 \ge r_k$ which forces $r_k=0$.

Now assume that $k$ is not minimal and let $m=M(k)$.  Since $m<_P k$ we have $\max\rho_m < r_k$.  On the other hand, $m+1\nless_P r_k$ so that $\max\rho_{m+1} \ge r_k$.
Using the previous inequalities and equation~\eqref{k,k+1} gives $\max\rho_m < r_k \le \max\rho_m+ 1$.  This implies  $r_k =\max\rho_m+ 1$ as desired.
\eprf

Now we defined a map $\rp^{-1}:\Av_n(P_3)\ra\RGF_n$ by letting $\rp^{-1}(P)$ be the sequence $\rho$ defined by equation~\eqref{rp-1}.
Again, we need to worry about whether the map is well defined.
\begin{lem}
If $P\in\Av_n(P_3)$ then $\rho=\rp^{-1}(P)\in\RGF_n$.
\end{lem}
\bprf
We claim  that since $P$ is factorial we must have that $1$ is a minimal element.  Suppose this were not the case.   Then the set of elements below $1$ would not form an inverval of the form $[m]$ for some $m>2$ since all such intervals contain $1$, which contradicts~\eqref{MkInt}.  But $1$ being minimal means that $\rho_1=0$ as needed for an RGF.

We now prove that (r2) in the definition of an RGF holds.  If $k$ is minimal then $r_k=0$ and the inequality is trivial.  If $k$ is not minimal then let $m=M(k)$.  Since $P$ is factorial we have $m\le k-1$.
It follows that $\max\rho_m\le\max\rho_{k-1}$.  Combining this with~\eqref{rp-1} gives
$$
r_k = \max\rho_m+1 \le\max\rho_{k-1} +1
$$
which is what we wished to prove.
\eprf

We can now show that the maps $\rp$ and $\rp^{-1}$ are indeed inverses.  We compose functions from right to left.
\bth
The map $\rp:\RGF_n\ra\Av_n(P_3)$ is a bijection with inverse $\rp^{-1}$.
Thus
$$
\#\Av_n(P_3) = B_n.
$$
\eth
\bprf
By~\eqref{RGFB}, the cardinality result will follow once we know that $\rp$ is  a bijection.

Since equation~\eqref{rp-1} was shown to hold for $\rp$ and then used to define $\rp^{-1}$, we automatically have that $\rp^{-1}\circ\rp$ is the identity map.  For the other composition,
consider $P\in\Av_n(P_3)$ and let $\rho=\rp^{-1}(P)$ and $P'=\rp(\rho)$.  By induction on $n$ we can assume that $P$ equals $P'$ when both are restricted to $[n-1]$.  There are now two cases depending on the element $n$.  Note that since $P$ and $P'$ are both naturally labeled, it suffices to consider the elements below $n$ in the two posets.

Suppose first that $n$ is minimal in $P$.  By the definition of $\rp^{-1}$ we have $r_n=0$.  But then there is no $i$ with $\max \rho_i < r_n$.  So $n$ is also minimal in $P'$ and we are done.

If $n$ is not minimal in $P$ then let $m=M(n)$ computed in $P$.  We need to show that $m$ is the element of maximum label below $n$ in $P'$.
First note that by~\eqref{rp-1} we have 
$$
r_n = \max \rho_m +1.
$$
So, in particular, $\max\rho_m < r_n$ and the definition of $\rp$ implies $m<_{P'} n$.  Thus we will be done if we can prove that $m+1\nless_{P'} n$.
If, instead, we have $m+1<_{P'} n$ then this is equivalent to $\max\rho_{m+1} < r_n$.  Combining this with the previous displayed equation gives
$$
\max\rho_m \le \max\rho_{m+1} < r_n=\max \rho_m +1.
$$
It follows that 
$$
\max\rho_{m+1} = \max\rho_m.
$$
Now recall that $P$ has no special $P_3$.  But $m<_P n$ and $m+1\nless_P n$ so we must have $m+1>_P m$ otherwise $m$, $m+1$, and $n$ form a special $P_3$.  And since $m,m+1$ are consecutive integers this forces $m = M(m+1)$ computed in $P$.  Now~\eqref{rp-1} implies $r_{m+1} = \max\rho_m + 1$.  This in turn yields $\max\rho_{m+1} = \max\rho_m + 1$ which directly contradicts the previous displayed equation.
\eprf

We will now make a connection between RGFs and pattern avoidance.
Claesson~\cite{cla:gpa} showed that the number of permutations in $S_n$ avoiding any of the following vincular patterns is the  Bell number, $B_n$:
$$
1|23, 12|3, 1|32, 21|3, 23|1, 3|12, 3|21, 32|1
$$
These are trivially Wilf equivalent to the bivincular patterns
$$
\ol{1}23, 1\ol{2}3, \ol{1}32, \ol{2}13, 23\ol{1}, 31\ol{2}, 3\ol{2}1, 32\ol{1}.
$$
We will translate a version of one of his bijections from set partitions into the language of RGFs.  Given $\rho=r_1 r_2 \ldots r_n\in\RGF_n$ we form a permutation $\pi=p_1 p_2\ldots p_n \in S_n$ as follows.
Suppose the indices  
 for which  $r_i=0$ are
\beq
\label{i's}
1=i_1 < i_2 <\ldots < i_k.
\eeq
Then $\pi$ will start with the initial segment
$$
\pi_k = i_k i_{k-1}\ldots i_1.
$$
Now let the indices  
 for which $r_j=1$ be
\beq
\label{j's}
j_1 < j_2 <\ldots < j_l.
\eeq
Continue building $\pi$ so that
$$
\pi_{k+l} = \pi_k j_l j_{l-1} \ldots j_1.
$$
The rest of $\pi$ is built 
in a similar manner.
\bpr[\cite{cla:gpa}]
The map $\rho\mapsto\pi$ is a bijection $\RGF_n\ra S_n(23|1)$.\hqed
\epr

There is also a standard way to associate with any RGF a strictly upper-triangular, $0$-$1$ matrix.  This is usually expressed in the language of set partitions and rook placements on a triangular Ferrers board.  We mention it here in translation.  Given $\rho\in\RGF_n$ we will associate an $n\times n$ matrix $M$ with rows and columns indexed by $[n]$.  Start with $M$ being the zero matrix.  Now consider the indices in~\eqref{i's} and for each adjacent pair $(i,j)=(i_a,i_{a+1})$ set $M_{i,j}=1$ giving $k-1$ nonzero entries.  Continue with the sequence~\eqref{j's}, letting $M_{i,j}=1$ for each pair
$(i,j)=(j_b,j_{b+1})$, and so forth.
\bpr
The map $\rho\mapsto M$ is a bijection between $\RGF_n$ and $n\times n$ strictly upper-triangular, $0$-$1$ matrices. \hqed
\epr

One can also associate with each RGF a perfect matching by restricting the Claesson-Linusson bijection{\bf~\cite{CL:n!m}} to the inversion sequences which are RGFs.  
However, the description of the image seems a bit contrived and so we will leave the details to the reader.

\section{Comments and open questions}
\label{coq}

We collect here some of the problems suggested by our work.

\subsection{More on matrices}

Notice that Theorem~\ref{MtxThm} only applies for $d\ge1$.  But a similar result holds for $d=0$ and the set of ascent sequences $\A_n$.  Specifically, factor an ascent sequence $\al=\de_0\de_1\ldots\de_m$ by breaking at every ascent.  So the $\de_i$ are maximal weakly decreasing factors of $\al$.  Now use~\eqref{M_k} to construct a matrix $\mx(\al)$.  Note that since the $\de_i$ weakly decrease, it is possible to add the same $E_{i,j}$ multiple times so that the entries  of the final matrix can be larger than one.  In fact, it is not hard to see that the image of the map is the set $\Mtx_n$ of all matrices $M$ who entries sum to $n$ and have the following properties.
\ben
\item[(Ma)]  $M$ is upper triangular with nonnegative integer entries.
\item[(Mb)]  There are no zero columns in $M$, and for all column indicies $j$ 
$$
\rmax c_j > \rmin c_{j-1}.
$$
\een
Thus, we have the following result
\bth
For $n\ge0$ the map $\mx:\A_n\ra\Mtx_n$ is a bijection.  Consequently

\vs{10pt}

\eqed{\#\A_n =\#\Mtx_n}
\eth

Dukes and Parviainen~\cite{DP:asu} have shown that ascent sequences are in bijection with another set of upper triangular matrices. Let $\Mtx'_n$ be the set of all matrices $M$ whose entries sum to $n$ and which satisfy (Ma) and 
\ben
\item[(Mc)]  There are no zero rows or columns in $M$.
\een
\bth[\cite{DP:asu}]
There is a bijection $\mx':\A_n\ra\Mtx'_n$.\hqed
\eth

\begin{problem}
Find a direct bijection $\Mtx_n\ra\Mtx'_n$ without composing $\mx^{-1}$ and  $\mx'$.
\end{problem}

\subsection{More on permutations and posets}

The inequalities in Theorems~\ref{AvSi} and~\ref{AvP} can be strict for $d\ge2$.
This raises the question of how different $\#\dA_n$ can be from $\#\Av_n(\si_{d+3})$ and from $\#\Av_n(P_{d+3})$.  One could also ask if equality could be obtained by enlarging the set of patterns or posets to be avoided.  Specifically, we have the following questions.
\begin{question}
Fix $d\ge2$.
\ben
\item[(a)] Does $\lim_{n\ra\infty} \#\dA_n/\#\Av_n(\si_{d+3})$ exists and, if so, what is its value?
\item[(b)] Is there a set $\Si_d$ of bivincular patterns containing $\si_d$ such that
$\#\dA_n=\#\Av_n(\Si_{d+3})$ for $n\ge0$?
\item[(c)]  Answer the analogues of (a) and (b) for  posets.
\een
\end{question}

\subsection{More on rooted duplication trees}

Recall that Theorem~\ref{RDTThm} only holds for $d=2$.  So a natural question to ask is the following.
\begin{question}
For $d\ge3$, is there a generalization $\dRDT_n$ of the set $\RDT_n$ (rooted duplication trees with $n$ leaves) which is in bijection with $\dI_{n-1}$ ($d$-increasing $d$-ascent sequences of length $n-1$)?
\end{question}

We note that there is an obvious analogue of rooted duplication trees which model the duplication of a  factor in a gene sequence $d$ times rather than just twice.
This results in a set of rooted $d$-ary trees with labeled leaves.  Unfortunately, these tree are more numerous than  $d$-increasing $d$-ascent sequences.

\subsection{Generating functions}

\begin{table}
    \centering
    \begin{tabular}{c||ccccccccccc}
    $d\setm n$& 0 & 1 & 2 & 3 & 4 & 5 & 6 & 7 & 8 & 9 & 10 \\ \hline\hline
    0   & 1 & 1 & 2 & 5 & 15 & 53 & 217 & 1014 & 5335 & 31240 & 201608 \\ \hline
    1   & 1 & 1 & 2 & 6 & 23 & 106 & 567 & 3440 & 23286 & 173704 & 1414102 \\ \hline
    2   & 1 & 1 & 2 & 6 & 24 & 118 & 682 & 4506 & 33376 & 273200 & 2444274 \\ \hline
    3  & 1 & 1 & 2 & 6 & 24 & 120 & 714 & 4896 & 37854 & 324792 & 3055320 \\ \hline
    4  & 1 & 1 & 2 & 6 & 24 & 120 & 720 & 5016 & 39624 & 348840 & 3378192 \\ \hline
    5  & 1 & 1 & 2 & 6 & 24 & 120 & 720 & 5040 & 40200 & 358800 & 3534120 \\ \hline
    6  & 1 & 1 & 2 & 6 & 24 & 120 & 720 & 5040 & 40320 & 362160 & 3600720
    \end{tabular}
    \caption{The values of  $\#\dA_n$ for $d\le 6$ and $n\le 10$}
    \label{dAnTab}
\end{table}

Let $a_n=\#\A_n$, the number of ascent sequences of length $n$.  In~\cite{BMCDK:2fp} the generating function for $a_n$ was derived using the kernel method.  This series had been shown earlier by Zagier~\cite{zag:vis}
to count matchings with neither left nor right nestings.
\bth[\cite{BMCDK:2fp}]
We have

\vs{10pt}

\eqed{\sum_{n\ge0}  a_n t^n =    \sum_{n\ge0} \prod_{i=1}^n \left( 1 - (1-t)^i\right). }    
\eth

\begin{problem}
Find the generating function for $\#\dA_n$ for all $d\ge0$.  A table of these values for $d\le 6$ and $n\le 10$ is given in~\ref{dAnTab}  
\end{problem}

It is not hard to at least find a functional equation that might be of use to solve the previous problem.  Suppose $d=1$ so we are considering weak ascents.  Define
$$
W(u,v)=W(t;u,v) = \sum_{n\ge 0}\sum_{\al\in\wA_n} t^n u^{\wasc\al} v^{a_n}
$$
where $a_n$ is the last element of $\al$.  The next result can be easily derived by considering how the possible last elements of $\al$ affect the number of weak ascents.
\bpr
We have
$$
[v-1+t(u-1)] W(u,v) = [(1+t)(v-1)+tu(1-v^2)]- t W(u,1) + tuv^2 W(uv,1).\ \qed
$$
\epr

It would also be nice to have a generating function for $\di_n$, the number of $d$-increasing $d$-ascent sequences of length $n$, and not just the recursion in Theorem~\ref{dinRec}.  Note that for $d=0$ or $d=1$ this sequence is the all ones or Catalan number sequence, respectively, whose generating functions are well known.
\begin{problem}
Find the generating function for $\di_n=\#\dI_n$ for all $d\ge0$.  A table of these values for $d\le 6$ and $n\le 10$ is given in~\ref{dInTab}  
\end{problem}

Enumeration of lattice paths may be helpful with this problem.
Take $\al=a_1 a_2 \ldots a_n\in \dI_n$ and form the sequence
$\al' = a_1' a_2' \ldots a_n'$ where $a_k' = a_k + (k-1)(d-1)$ for $k\in[n]$.
From properties (I1) and (I2) one sees that 
\ben
\item[(I1')] $a_1'=0$, and
\item[(I2')] $a_k'\le a_{k+1}'\le kd$ for $k\in[n-1]$
\een
Letting $a_k'$ be the height of the $k$th east step, one obtains a lattice path using north and east steps and staying inside a cone.  It may also be advantageous to do a linear transformation so that the paths in question are restricted to the first quadrant, although now the two steps used will not be as simple as in the conical case.

\begin{table}
    \centering
    \begin{tabular}{c|| ccccccccccc}
    $d\setm n$  & 0 & 1 & 2 & 3 & 4 & 5 & 6 & 7 & 8 & 9 & 10 \\ \hline\hline
    0           & 1 & 1 & 1 & 1 & 1 & 1 & 1 & 1 & 1 & 1 & 1 \\ \hline
    1           & 1 & 1 & 2 & 5 & 14 & 42 & 132 & 429 & 1430 & 4862 & 16796 \\ \hline
    2           & 1 & 1 & 2 & 6 & 22 & 92 & 420 & 2042 & 10404 & 54954 & 298648 \\ \hline
    3           & 1 & 1 & 2 & 6 & 24 & 114 & 612 & 3600 & 22680 & 150732 & 1045440 \\ \hline
    4           & 1 & 1 & 2 & 6 & 24 & 120 & 696 & 4512 & 31920 & 242160 & 1942800 \\ \hline
    5           & 1 & 1 & 2 & 6 & 24 & 120 & 720 & 4920 & 37200 & 305280 & 2680800 \\ \hline
    6           & 1 & 1 & 2 & 6 & 24 & 120 & 720 & 5040 & 39600 & 341280 & 3175200
    \end{tabular}
    \caption{The values of  $\#\dI_n$ for $d\le 6$ and $n\le 10$}
    \label{dInTab}
\end{table}

\nocite{*}
\bibliographystyle{alpha}

\end{document}